\documentclass[review,hidelinks,onefignum,onetabnum]{siamart250211}

\usepackage{lipsum}
\usepackage{amsfonts}
\usepackage{graphicx}
\usepackage{epstopdf}
\usepackage{algorithmic}
\usepackage{tikz}
\usepackage{enumitem}
\usepackage{mathtools}
\usepackage{booktabs}
\usepackage{amsopn}

\ifpdf
  \DeclareGraphicsExtensions{.eps,.pdf,.png,.jpg}
\else
  \DeclareGraphicsExtensions{.eps}
\fi

\newsiamremark{remark}{Remark}
\newsiamremark{hypothesis}{Hypothesis}
\crefname{hypothesis}{Hypothesis}{Hypotheses}
\newsiamthm{claim}{Claim}
\newsiamremark{fact}{Fact}
\crefname{fact}{Fact}{Facts}

\headers{Practical Semialgebraic Optimization}{J. Nezvadovitz and D. Rosen}

\title{Practical Semialgebraic Optimization via Sparse Sample-Based Sum-of-Squares Programming\thanks{This manuscript will be submitted to the SIAM Journal on Optimization.
\funding{MIT-LL and the OUSD(R\&E) through Air Force Contract~FA8702-25-D-B002.}}}

\author{Jason Nezvadovitz\thanks{MIT Lincoln Laboratory, Lexington, MA
  (\email{nez.jason7@gmail.com}, \url{https://www.ll.mit.edu}).}
\and David M. Rosen\thanks{Departments of Mathematics, Electrical \& Computer Engineering, Northeastern University, Boston, MA
  (\email{d.rosen@northeastern.edu}, \url{https://neural.lab.northeastern.edu}).}
}

\let\savedreviewnotice\footercopyright
\renewcommand{\footercopyright}{\phantom{\savedreviewnotice}}

\ifpdf
\hypersetup{
  pdftitle={Practical Semialgebraic Optimization via Sparse Sample-Based Sum-of-Squares Programming},
  pdfauthor={J. Nezvadovitz and D. Rosen}
}
\fi

\begin{document}

\maketitle

\begin{abstract}
As the importance of semialgebraic optimization grows in engineering, so does the need for scalable sum-of-squares (SOS) methods.
Their convex relaxations can guide mixed-integer search and often deliver certified global optima outright.
Scaling traditional SOS typically requires hand-crafted bases and redundant constraints, whereas sample-based SOS numerically exploits feasible set geometry to reduce bases and build smaller, often tighter relaxations from black-box evaluations.
However, discrete variables can require exponentially many samples, and each sampled equation couples all cliques, destroying correlative sparsity.
In this work, we avoid enumeration and slack variables by sampling a tractable variety that contains the feasible set, with multipliers only for the omitted constraints.
To recover sparsity, we introduce \textit{separator polynomials} that split the sampled equations into local systems coupled only on clique overlaps.
This yields a novel, convergent hierarchy of semidefinite programs (SDPs) with block sizes and sample requirements governed by local interactions rather than total problem size.
Lastly, we show how to verify sample sufficiency and prune redundant equations numerically before solving.
In robust estimation and motion planning, we obtain certifiably small optimality gaps at scales far beyond those reported for traditional sparse SOS.
For example, we solve a truncated least-squares 3D registration problem with $1600$ candidate correspondences in under $49$ seconds using an off-the-shelf SDP solver on a laptop CPU.
\end{abstract}




\section{Introduction}
\label{sec:introduction}

Across engineering optimization, from structural design and power-grid operation to sensor fusion and motion planning, there is a growing recognition that many models' objective and constraints admit exact polynomial representations \cite{handa2025structural,ghaddar2015power,papalia2024slam,kang2024trajopt}.
Such optimization problems are \textit{semialgebraic}, a property that enables systematic construction of remarkably tight convex relaxations \cite{parrilo2012book}.
These relaxations can initialize local methods, provide bounds for global search, and in many cases deliver certified optima outright.
Moreover, their unified treatment of continuous and discrete variables offers an alternative avenue for mixed-integer problems when local heuristics are brittle and tree search is prohibitive \cite{yang2022outliers}.

The key enabling machinery is the \textit{Positivstellensätze}, a class of theorems that characterize when a polynomial is positive on a set specified by polynomial constraints.
As we will explain in \cref{sec:foundation}, such theorems recast tightening a bound on the objective as a \textit{sum-of-squares (SOS) program}, a convex optimization over the infinite-dimensional cone of SOS polynomials.
Then by parameterizing an inner approximation of this cone, one can obtain a \textit{relaxation hierarchy}, a sequence of finite-dimensional convex programs that tighten the bound monotonically \cite{parrilo2003hierarchy}.

Unfortunately, the general NP-hardness of semialgebraic optimization forces a tradeoff between a hierarchy's convergence and its scalability.
Those with guaranteed (and often finite) convergence, such as Lasserre’s hierarchy and BSOS \cite{lasserre2001hierarchy,weisser2018bsos}, yield combinatorially large semidefinite programs (SDPs), limiting their use to small or very sparse problems \cite{majumdar2020survey,marandi2018bsoseval}.
Scalable alternatives like RLT and DSOS/SDSOS \cite{sherali1992rlt,ahmadi2019sdsos} replace SDPs with linear or second-order cone programs, but they typically lack convergence guarantees and the tightness required for solution extraction \cite{zheng2019sdsoseval}.
In order to uniformly improve both convergence and scalability, additional structure or information must be leveraged.

One such lever is found in problem representation.
The polynomials defining the feasible set are not unique, and the chosen representation strongly influences the performance of a given relaxation hierarchy \cite{laurent2008book}.
\textit{Sample-based} SOS leverages a finite set of feasible points as a purely geometric representation.
Assuming the points are \textit{poised}, it yields valid relaxations at least as tight as those from the best algebraic representation \cite{cifuentes2017sampling}.
Remarkably, these relaxations are also lighter-weight because they eliminate multiplier variables and reduce the objective's degree to its effective value on the feasible set.
Additionally, the points enable computation of a minimal orthonormal basis for polynomials on the feasible set, further improving dimensionality and numerical conditioning.
Finally, because sample-based SOS only evaluates functions at these points, the objective and constraints can be given as black-boxes, avoiding large-scale polynomial expansions and making the software implementation easier to use.

\begin{figure}[t]
    \centering
    \includegraphics[width=\linewidth]{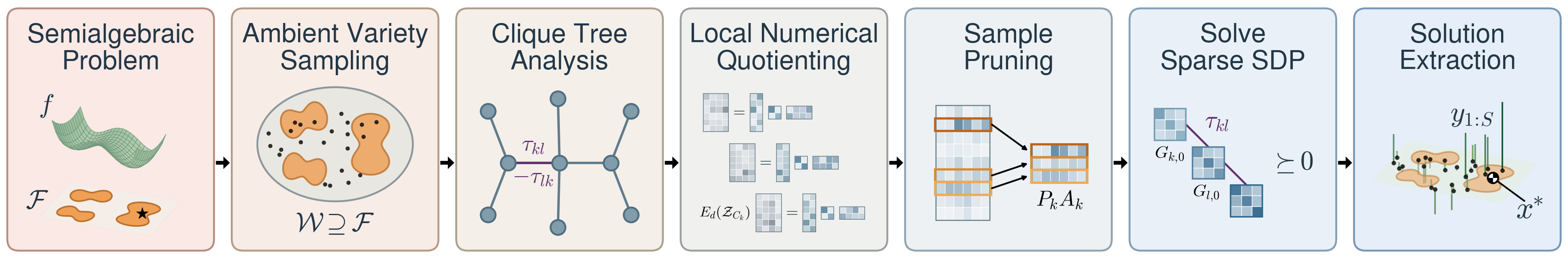}
    \vspace{-1em}
    \caption{
    Conceptual depiction of our methodology.
    The user selects a tractable variety $\mathcal{W}$ containing their problem's feasible set $\mathcal{F}$ and evaluates their objective and constraints at random samples on it.
    From there, correlative sparsity guides clique decomposition and the assignment of separator polynomials to selected clique overlaps.
    SVDs of local monomial evaluation matrices construct efficient certificate and separator bases by quotienting out redundant directions.
    Then the sampled certificate equations are assembled into sparsely coupled local systems, with rank checks verifying sample sufficiency before excess equations are pruned.
    Finally, the resulting sparse SDP is solved, and the solution's dual sample moments reveal the original problem's global optima.
    }
    \label{fig:overview}
    \vspace{-0.25em}
\end{figure}

To obtain a poised collection of feasible points, prior work assumes that the feasible set is an \textit{irreducible algebraic variety} --- i.e., defined solely by polynomial equalities and not equal to any finite union of other varieties.
Under these conditions, poisedness is readily achieved by random sampling from any distribution fully supported on the feasible set.
In theory, this assumption is without loss of generality: inequalities can be converted to equalities by adding squared slack variables, and reducible varieties can be handled componentwise.

In practice, however, the limitations are substantial.
For problems with many inequalities (e.g.\ control-limited motion planning), slack variables bloat the ambient dimension, largely negating any dimensionality reduction benefits.
And for problems with many discrete variables (e.g.\ tasking and robust estimation), solving a relaxation for every possible assignment is entirely intractable.
Moreover, sample-based SOS appears at odds with sparsity because sample evaluation mixes the coefficient structure used by traditional sparse methods, leaving the sampled system globally coupled even when all interactions in the original problem are local.
Consequently, despite its favorable numerical properties, sample-based SOS has been limited in practice to small dense problems with few discrete variables or inequality constraints.

In this work, we extend sample-based SOS to general sparse semialgebraic problems without introducing slack variables or enumerating irreducible components.
Our first insight is that convergence and scalability benefits can be derived from poised samples of a variety \textit{containing} the feasible set.
This enables us to use samples that satisfy a subset of the constraints (e.g.\ omitting inequalities and integrality) while still exploiting the structure of the variety defined by the remaining constraints.
Then we show how correlative sparsity can be leveraged through \textit{separator polynomials} to decompose both the SOS basis construction and the resulting convex program, all while preserving black-box problem specification.
Finally, we provide a simple method for pruning samples to remove linearly dependent constraints and verify poisedness before the solve.
An overview is shown in \cref{fig:overview}.

In effect, our approach unifies the numerical advantages of sample-based SOS with the direct applicability of traditional SOS.
We evaluate it on robust estimation and motion planning problems that were previously out of reach for sample-based SOS due to numerous discrete variables and inequality constraints.
By leveraging both sparsity and samples, our relaxations solve these problems at scales far beyond those demonstrated by the leading traditional sparse SOS approaches.
These results indicate a new state of the art for semialgebraic optimization, reducing reliance on hand-crafted formulations in favor of practical black-box evaluations.


\subsection{Related Work}
\label{sec:related_work}

Sample-based SOS originated with \cite{lofberg2004sampling} for unconstrained polynomial optimization, revealing improved numerical properties and a convenient black-box interface.
Later, \cite{cifuentes2017sampling} extended the idea to equality-constrained problems by sampling every irreducible component of the feasible variety, substantially broadening its scope but still excluding mixed-integer and inequality-heavy settings.
Moreover, while they suggested that standard sparsity-exploiting techniques could be applied, sample evaluation masks the structure that those methods rely on, leaving the sampled certificate equations globally coupled.
These limitations have curtailed practical adoption: to date, sample-based SOS has appeared primarily on small problems with irreducible algebraic feasible sets (e.g., \cite{shen2020lyap,zhang2023cbf}).
Meanwhile, the growing body of applications with critical sparsity and combinatorial reducibility continues to rely on bespoke algebraic manipulations for traditional SOS (e.g., \cite{yang2022outliers,kang2024trajopt}).
This gap between the promise of sample-based SOS and its current applicability motivates our work.

There have been a few other approaches to generalizing \cite{cifuentes2017sampling} that differ from ours.
Autotight \cite{dumbgen2024autotight} uses samples to compute a strong algebraic description of the feasible set, then propagates this description to related problem instances through a template-and-prune procedure.
This avoids resampling for each new instance but forfeits the dimensionality reduction of true sample-based SOS and still presumes an irreducible feasible set in practice.
Kernel-based SOS \cite{rudi2025kernelsos} generalizes sample-based SOS to arbitrary optimization by replacing Positivstellensatz certificates with a representer-theorem that yields an SDP in a non-polynomially lifted space.
However, convergence is not guaranteed and the resulting SDP is dense with cone size coupled to the sample count rather than intrinsic dimension, necessitating a ``restart strategy'' to guide the restricted sample-budget \cite{groudiev2025kernelsos}.
These directions reveal interesting options for black-box solvers but do not deliver the combination of scalability and convergence guarantees that we target.


\subsection{Contribution Outline}
\label{sec:contribution_outline}

We begin in \cref{sec:foundation} with a self-contained treatment of the theory underlying our work.
To keep the exposition broadly accessible, we assume no prior knowledge of semialgebraic geometry and provide illustrative examples throughout.
We also focus on real-valued decision variables for simplicity, though the complex case is no harder in principle.
Then, building on this foundation, we present the following novel extensions to sample-based SOS:

\begin{itemize}[itemsep=0.75em, topsep=0.75em, parsep=0pt, partopsep=0pt]
    \item In \cref{sec:circumventing_enumeration}, we show how to apply sample-based SOS to reducible problems by sampling from a larger irreducible variety that contains the feasible set. This enables solving mixed-integer problems without brute-force enumeration of discrete assignments.
    \item In \cref{sec:handling_inequalities}, we incorporate polynomial inequality constraints directly into the sample-based framework, avoiding slack-variable bloat while retaining convergence guarantees.
    \item In \cref{sec:leveraging_sparsity}, we exploit correlative sparsity in a novel way that decomposes both the monomial evaluation matrix and the entire convex program, yielding block-sparse SOS relaxations from a fully sample-based interface.
    \item In \cref{sec:selecting_samples}, we propose a method for reducing the sample set to eliminate linearly redundant constraints and verify poisedness before the solve.
\end{itemize}

Finally, in \cref{sec:results} we deliver our methodology as an efficient, publicly available \textbf{\href{https://github.com/Rhiza-PolyOpt}{codebase}} and evaluate it on large-scale robust estimation and motion planning problems.
With far less problem-specific algebraic manipulation, our procedure automatically produces SDPs that are substantially smaller than existing hand-tailored formulations.
Even on a laptop CPU, this reduction translates into \textbf{orders of magnitude faster solves at substantially larger scales} than previously reported for both problem classes.
We aim for this work to be of immediate use to optimization practitioners, enabling large-scale certifiably-global semialgebraic optimization from \textbf{black-box evaluations rather than per-problem algebraic engineering}.


\subsection{Notation}
\label{sec:notation}

For $i,j \in \mathbb{N} := \{0,1,\ldots\}$, we use $a_{i:j}$ to denote the ordered tuple $(a_i,\ldots,a_j)$, considered empty if $j<i$.
Relations written with this shorthand apply elementwise.
For example, $x_{1:n} \in \mathbb{R}$ means that each $x_i$ is real, and $g_{1:m}(x) \geq 0$ requires $g_i(x) \geq 0$ for all $i \in \{1,\ldots,m\}$.
With an index set $C \subseteq \{1,\ldots,n\}$ of cardinality $|C|$, we write $x_C$ for the corresponding subvector of $x \in \mathbb{R}^n$ and $\mathbb{R}[x_C]$ for the ring of real polynomials in those variables.
For $\alpha \in \mathbb{N}^n$, we let $x^\alpha := \prod_{i=1}^n x_i^{\alpha_i}$ be a monomial of total degree $|\alpha| := \sum_{i=1}^n \alpha_i$.
In the space $\mathbb{S}^N$ of real symmetric $N \times N$ matrices, $G \succeq 0$ denotes membership in the positive semidefinite cone $\mathbb{S}^N_+$.
The length-$S$ vector of ones is $\mathbf{1}_S$.
Finally, $\operatorname{svec}$ vectorizes symmetric matrices so as to preserve the Frobenius inner product: $\operatorname{svec}(A)^\top\!\operatorname{svec}(B) = \langle A,B \rangle := \operatorname{tr}(A^\top\!B)$.


\section{Foundation}
\label{sec:foundation}

\begin{figure}[htbp]
    \centering
    \includegraphics[width=0.7\linewidth]{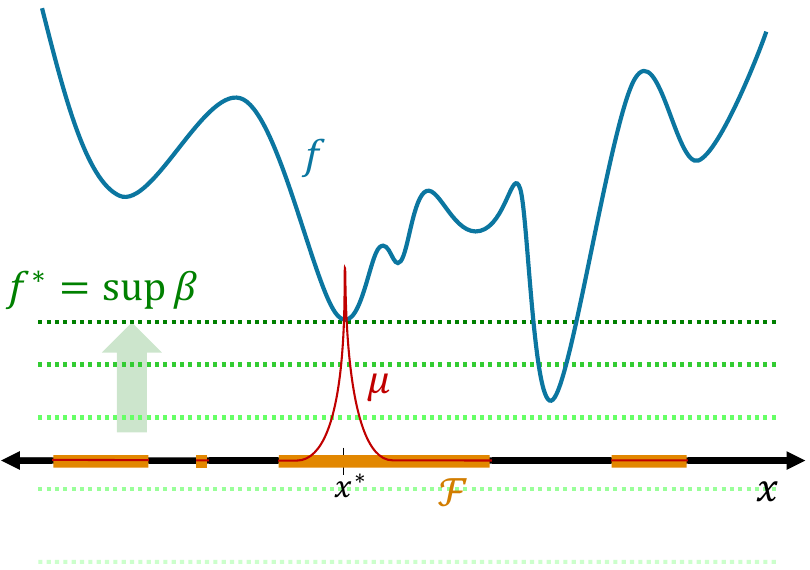}
    \caption{
    Depiction of a minimization problem $(f,\mathcal{F})$, its GLB $\beta$, and dual measure $\mu$, in 1D $(\mathcal{F} \subset \mathbb{R})$ with unique minimizer $x^*$.
    The objective function $f$ shown in blue is minimized over the feasible set $\mathcal{F}$ highlighted in orange.
    The minimum value $f^*$ is the supremum of lower bounds $\beta$, depicted by dashed green lines ascending toward $f^*$.
    Note that $\beta$ can exceed $f$ at any $x \notin \mathcal{F}$.
    Each point in $\mathcal{F}$ corresponds to an inequality constraint on $\beta$ and thus a nonnegative dual value.
    As $\beta$ approaches $f^*$, the red distribution of dual values approaches a Dirac delta at $x^*$ for the sole active constraint $\beta \leq f(x^*)$.
    For multiple global minima, $\mu$ may concentrate on any subset of them.
    }
    \label{fig:glb}
\end{figure}

Let $\mathcal{F}$ be the feasible set and $f : \mathcal{F} \to \mathbb{R}$ be the objective function of the general optimization problem,
\begin{equation}
    f^* \,:=\, \inf_{x \in \mathcal{F}}\, f(x)
    \label{eq:opt}
\end{equation}
A lower bound is a value $\beta \in \mathbb{R}$ that does not exceed the objective $f(x)$ at any feasible point $x \in \mathcal{F}$.
Evidently, the optimal value $f^*$ is the greatest lower bound (GLB),
\begin{equation}
    f^* \,=\, \sup_{\beta \in \mathbb{R}}\, \beta
    \quad \text{s.t.} \quad
    \beta \leq f(x) \;\ \forall x \in \mathcal{F}
    \label{eq:glb}
\end{equation}
Through this GLB lens, a general optimization problem becomes a linear program (LP) in the single variable $\beta$, with one constraint for each point in $\mathcal{F}$.
As illustrated in \cref{fig:glb}, any such constraint active at optimality corresponds to a global minimizer of the original problem \cref{eq:opt}.
This suggests examining the dual LP, whose variables place weight on the active constraints and thereby indicate minimizers.

The dual variable is a measure $\mu$ supported on $\mathcal{F}$, which can be interpreted as assigning a nonnegative Lagrange multiplier to each inequality constraint $\beta \leq f(x)$ indexed by $x \in \mathcal{F}$.
The Lagrangian is
\begin{align*}
    \mathcal{L}(\beta,\mu) \,:=\,&\; \textcolor{gray}{\overbrace{\normalcolor \beta}^{\text{objective}}} +\; \int_\mathcal{F} \textcolor{gray}{\overbrace{\normalcolor \big(f(x) - \beta\big)}^{\text{constraints}}} \textcolor{gray}{\overbrace{\normalcolor d\mu(x)}^{\text{multipliers}}} \\
    =\,&\; \int_\mathcal{F} f(x)\, d\mu(x) \,+\, \big(1 - \mu(\mathcal{F})\big) \beta
\end{align*}
and consequently the dual objective is
\begin{equation*}
    \sup_{\beta \in \mathbb{R}}\, \mathcal{L}(\beta,\mu) \,=\,
    \begin{cases}
        \int_\mathcal{F} f(x)\, d\mu(x) & \text{if }\ \mu(\mathcal{F}) = 1 \\
        \infty & \text{otherwise }
    \end{cases}
\end{equation*}
So dual feasibility\footnote{
    Normalization of $\mu$ can be equivalently understood as the KKT stationarity condition $\frac{\partial \mathcal{L}}{\partial \beta} = 0$.
} amounts to the normalization condition $\mu(\mathcal{F}) = 1$, i.e.\ it must be a probability measure.
Thus we see the dual LP as minimizing the expected objective value over all probability measures supported on the feasible set $\mathcal{P}(\mathcal{F})$,
\begin{equation}
    f^* \,=\, \inf_{\mu \in \mathcal{P}(\mathcal{F})}\, \underset{x \sim \mu}{\mathrm{E}}\big[f(x)\big]
    \label{eq:gmp}
\end{equation}

To make lower objective values more likely, any optimal measure $\mu^*$ concentrates its mass on the minimizers of \cref{eq:opt}.
Recovering a solution therefore reduces to locating that mass.
For example, if the minimum of \cref{eq:opt} is attained uniquely at $x^*$, then $\mu^*$ is the Dirac measure $\delta_{x^*}$ and we can extract the minimizer from its mean: $x^* \,=\, \underset{x \sim \mu^*}{\mathrm{E}}\big[x\big]$.

\begin{figure}[htbp]
    \centering
    \includegraphics[width=0.75\linewidth]{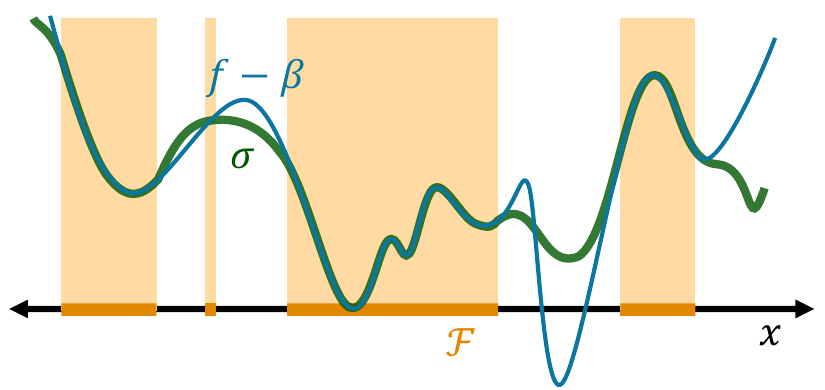}
    \caption{
    A candidate lower bound $\beta$ is valid exactly when the function $f-\beta$ is nonnegative on $\mathcal{F}$.
    Rather than check this condition pointwise, we seek a globally nonnegative function $\sigma$, shown here in green, that matches $f-\beta$ on $\mathcal{F}$.
    This certificate guarantees $f-\beta \geq 0$ on $\mathcal{F}$ even if it is negative outside $\mathcal{F}$.
    Maximizing $\beta$ pulls $f-\beta$ downward until it first touches zero on $\mathcal{F}$ at a global minimizer.
    Increasing $\beta$ any further makes $f-\beta$ negative somewhere on $\mathcal{F}$, so $\sigma$ ceases to exist.
    }
    \label{fig:cert}
\end{figure}

Thus \cref{eq:glb} is an LP whose primal and dual solutions provide the globally optimal value and decision variables for the general (nonconvex) problem \cref{eq:opt}.
However, in continuous and mixed-integer settings where $|\mathcal{F}|$ is infinite, how can we handle the infinite number of constraints?
The key idea is to view them collectively as a single functional constraint: $f - \beta \,\text{ is nonnegative on } \mathcal{F}$.
As illustrated in \cref{fig:cert}, this constraint is satisfied if we can find a globally nonnegative function $\sigma$ that matches $f - \beta$ on $\mathcal{F}$.
We call $\sigma$ a \textit{certificate} for the validity of $\beta$ as a lower bound.

In this sense, optimization \cref{eq:opt} is subsumed by the broader problem of finding nonnegativity certificates.
The set of nonnegative functions is a convex cone, but intractable to search in full.
Fortunately, there are broad problem classes $(f,\mathcal{F})$ for which theorems known as \textit{Positivstellensätze} allow us to restrict attention to structured subcones, reducing the semi-infinite LP to finite-dimensional conic programs.
This work is focused on one such class: semialgebraic optimization.


\subsection{Semialgebraic Optimization}
\label{sec:semialgebraic_optimization}

A \textit{semialgebraic optimization} problem is one whose epigraph
is a semialgebraic set, i.e., definable without quantifiers over the ordered real field.
More concretely, any such problem can be formulated through auxiliary variables as a \textit{polynomial program} (POP):
minimizing a real polynomial over a set specified by finitely many polynomial constraints.
That is, for some $n, m, \ell \in \mathbb{N}$ and polynomials $g_{1:m},\, h_{1:\ell} \,\in\, \mathbb{R}[x_{1:n}]$, we consider problems of form \cref{eq:opt} where the objective function and feasible set are given by
\begin{equation}
\begin{gathered} \label{eq:pop}
    f \,\in\, \mathbb{R}[x_{1:n}] \\[5pt]
    \mathcal{F} \,=\, \Big\{x \in \mathbb{R}^n \;\,\Big\vert\;\, g_{1:m}(x) \geq 0 \,,\ h_{1:\ell}(x) = 0\Big\}
\end{gathered}
\end{equation}

Many problems in engineering already have this form, including least-squares over rigid-body poses, similarity transforms, and complementarity constraints.
Many more can be converted to it through simple \textit{polynomial lifts}.
For example, a rational expression $p(x)/q(x)$ can be represented by an auxiliary variable $z$ satisfying the polynomial constraint $q(x)z - p(x) = 0$.
Likewise, a trigonometric expression such as $\sin(x)\cos^3(x)$ can be replaced by $sc^3$ together with $s^2 + c^2 - 1 = 0$.
Even discrete variables can be encoded as polynomial roots, e.g.\ $x(x-1)=0$ for $x \in \{0,1\}$, which in turn enables conditional logic such as \texttt{[if x then y=p(z)]} to be expressed as $x\big(y-p(z)\big)=0$.
Similar encodings allow \mbox{piecewise relations, permutations, and more.}

The equality constraints in \cref{eq:pop} generally control the effective dimension of the problem.
A subset of $\mathbb{R}^n$ determined by finitely many polynomial equalities is called a \textit{real algebraic variety},\footnote{
    In keeping with the optimization literature, we use this term in the affine sense and do not assume irreducibility unless stated otherwise.
}
which we abbreviate here as just ``variety.''
Accordingly, for the developments that follow, it is helpful to distinguish the variety $\mathcal{V} \supseteq \mathcal{F}$ defined by the equality constraints alone,
\begin{equation*}
    \mathcal{V} \,:=\, \Big\{x \in \mathbb{R}^n \;\,\Big\vert\;\, h_{1:\ell}(x) = 0\Big\} 
\end{equation*}

From a purely representational viewpoint, the distinction between equality and inequality constraints is not fundamental.
Any equality $p(x) = 0$ can be expressed as the pair of inequalities $p(x) \geq 0$ and $-p(x) \geq 0$, while any inequality $q(x) \geq 0$ can be converted to an equality by introducing a slack variable $s$ with $q(x)-s^2 = 0$.
For ease of explanation and to match the original sample-based SOS literature, we adopt the equalities-only standpoint for the next few sections by assuming that $x$ includes any necessary slack variables.
Equivalently, we temporarily restrict attention to the case $m=0$ for which $\mathcal{F}=\mathcal{V}$.
Although this is technically without loss of generality, we will show in \cref{sec:handling_inequalities} that preserving explicit inequality constraints is nevertheless important computationally.

Finally, it is crucial to recognize that the polynomials $h_{1:\ell}$ are a particular \textit{algebraic description} of $\mathcal{V}$, and such descriptions are not unique.
$\mathcal{V}$ exists as a geometric object independent of the specific equations supplied with the problem.
For example, the unit circle can be described equivalently by $x^2+y^2-1$ and $x^4+x^2y^2+y^2-1$.
This suggests shifting attention from $h_{1:\ell}$ to the set of all polynomial equalities implied by the variety $\mathcal{V}$ itself: its \textit{vanishing ideal}\footnote{
    An ideal of a commutative ring $\mathcal{R}$ is an additive subgroup of $\mathcal{R}$ closed under multiplication by any element of $\mathcal{R}$.
    Indeed for $\mathcal{I}(\mathcal{V}) \subseteq \mathbb{R}[x_{1:n}]$, if $p_1$ and $p_2$ vanish on $\mathcal{V}$, then so does $q_1p_1 + q_2p_2$ for any $q_1, q_2 \in \mathbb{R}[x_{1:n}]$.
    While vanishing ideals can be defined in any function ring, polynomial ideals are especially useful because they are finitely generated.
}\begin{equation*}
    \mathcal{I}(\mathcal{V}) \,:=\, \Big\{p \in \mathbb{R}[x_{1:n}] \;\,\Big\vert\;\, p(x) = 0 \;\ \forall x \in \mathcal{V}\Big\}
\end{equation*}
The vanishing ideal captures when two polynomials are indistinguishable on $\mathcal{V}$,
\begin{equation*}
    p(x) = q(x) \;\ \forall x \in \mathcal{V}
    \quad \iff \quad
    p - q \,\in\, \mathcal{I}(\mathcal{V})
\end{equation*}
Thus, in the equalities-only ($\mathcal{F} = \mathcal{V}$) setting, our problem
\begin{equation}
    f^* \,=\, \inf_{x \in \mathcal{V}}\, f(x)
    \label{eq:vopt}
\end{equation}
is naturally understood \textit{modulo} $\mathcal{I}(\mathcal{V})$: we can add any element of $\mathcal{I}(\mathcal{V})$ to the objective without changing the problem.
This will be useful when we compute nonnegativity certificates, as membership in polynomial ideals can be checked in finite time.


\subsection{Sum-of-Squares Programming}
\label{sec:sum_of_squares_programming}

Returning to the GLB view, a lower bound $\beta$ for \cref{eq:vopt} is valid when $f-\beta$ is nonnegative on the variety $\mathcal{V}$, which we enforce by jointly finding a globally nonnegative function $\sigma$ that matches $f-\beta$ on $\mathcal{V}$,
\begin{equation}
    f^* \,= \sup_{\substack{\beta \,\in\, \mathbb{R} \\ \sigma \,:\, \mathbb{R}^n \to \mathbb{R}_+}}\, \beta
    \quad \text{s.t.} \quad
    f(x) - \beta = \sigma(x) \;\ \forall x \in \mathcal{V}
    \label{eq:vglb}
\end{equation}
This form reveals two distinct challenges:
The first is tractably parameterizing global nonnegativity for $\sigma$.
The second is imposing agreement on $\mathcal{V}$.

For nonnegativity, it is natural to express $\sigma$ as a \textit{sum-of-squares} (SOS),
\begin{equation*}
    \sigma \,=\, \sum_i q_i^2\,, \quad q_i : \mathbb{R}^n \to \mathbb{R}
\end{equation*}
To search this representation, we choose a function basis $\varphi : \mathbb{R}^n \to \mathbb{R}^N$ and restrict each $q_i$ to lie in the $\operatorname{span}(\varphi_{1:N})$.
Using $c_i \in \mathbb{R}^N$ to write $q_i = c_i^\top \varphi$ we get
\begin{equation*}
\begin{aligned}
    \sigma(x) \,=\,& \textstyle\sum_i q_i(x)^2 \\[2pt]
    =\,& \textstyle\sum_i \big(c_i^\top \varphi(x)\big)^2 \\[2pt]
    =\,& \varphi(x)^\top \big(\textstyle\sum_i c_i c_i^\top\big)\, \varphi(x) \\[2pt]
    =\,& \varphi(x)^\top G\, \varphi(x)
\end{aligned}
\end{equation*}
where the \textit{Gramian} $G := \sum_i c_i c_i^\top$ is a positive semidefinite (PSD) matrix.
Conversely, every PSD matrix admits such a decomposition.
Thus, for a fixed basis $\varphi$, searching over SOS functions amounts to searching over PSD matrices $G \in \mathbb{S}^N_+$.

For agreement on $\mathcal{V}$, polynomial structure becomes key.
If we choose $\varphi$ to be polynomial, then $\sigma = \varphi^\top G\,\varphi$ is polynomial as well, and agreement on $\mathcal{V}$ becomes the membership $f - \beta - \sigma \in \mathcal{I}(\mathcal{V})$.
We will show later that this ideal membership condition can be expressed as a finite linear system in $\beta$ and $G$.
Therefore, restricting $\sigma$ to SOS polynomials reduces \cref{eq:vglb} to a semidefinite program (SDP),
\begin{equation}
    f^*_\varphi \,:= \sup_{\substack{\beta \,\in\, \mathbb{R} \\ G \,\in\, \mathbb{S}^N_+}}\, \beta
    \quad \text{s.t.} \quad
    \underset{\rule{0pt}{2ex}\text{\textcolor{gray}{(will be linear in $\beta$ and $G$)}}}{f - \beta - \varphi^\top G\,\varphi \;\in\; \mathcal{I}(\mathcal{V})}
    \label{eq:vsos}
\end{equation}

Every SOS polynomial is nonnegative, so any feasible certificate $\sigma = \varphi^\top G\,\varphi$ in \cref{eq:vsos} also certifies a valid lower bound in \cref{eq:vglb}.
Therefore $f^*_\varphi \leq f^*$: the SOS restriction may sacrifice bound tightness, but not validity.
Enlarging the basis span preserves all existing certificates while admitting new ones, so the bound cannot worsen.
An exhaustive nesting of bases thereby yields a sequence of monotonically tightened lower bounds, which is called an \textit{SOS relaxation hierarchy}.
\begin{equation*}
    \operatorname{span}(\varphi_{1:N}) \subset \operatorname{span}(\varphi'_{1:N'})
    \quad \implies \quad
    f^*_\varphi \leq f^*_{\varphi'} \leq f^*
\end{equation*}

The critical question then is whether such hierarchies converge all the way to $f^*$, or whether there is a gap.
This question can be answered by Positivstellensatz theorems, which relate positivity on constrained sets to algebraic certificates.
In~particular, we use the following result adapted from \cite[Corollary~3]{schmudgen1991positiv}.
\begin{theorem}[Positivstellensatz for Compact Varieties]
    \label{thm:vpositiv}
    Let $\mathcal{V}$ be a compact real algebraic variety.
    If a real polynomial $p$ is strictly positive on $\mathcal{V}$, then there exists an SOS polynomial $\sigma$ such that $p - \sigma \,\in\, \mathcal{I}(\mathcal{V})$.
\end{theorem}

The compactness requirement is natural in our setting: engineering decision variables typically have physical limits, so we assume throughout that $\mathcal{F}$ is compact.
Thus, in the current equalities-only setting where $\mathcal{F} = \mathcal{V}$, our variety $\mathcal{V}$ is compact and \cref{thm:vpositiv} applies.
For any $\beta < f^*$, the polynomial $p = f - \beta$ is strictly positive on $\mathcal{V}$, so the theorem guarantees existence of an SOS polynomial certificate, thereby ensuring asymptotic convergence from below ($f^*_\varphi \uparrow f^*$).

Even though \cref{thm:vpositiv} only guarantees certificates for $\beta < f^*$, finite convergence to $\beta = f^*$ occurs often in practice.
When it does, the certificate satisfies $\varphi^\top G^*\,\varphi = f - f^*$ on $\mathcal{V}$ and touches zero on $\mathcal{V}$ exactly at the global minimizers, just as in \cref{fig:cert}.
Since $G^* \succeq 0$, this reduces solution extraction to a nullspace problem,
\begin{equation*}
    x^* \in\, \operatorname*{arg\,min}_{x \,\in\, \mathcal{V}}\, f(x)
    \quad \iff \quad
    G^* \varphi(x^*) = 0\,,\;\; x^* \in \mathcal{V}
\end{equation*}
Furthermore, the dual variable $M \in \mathbb{S}^N_+$ paired with the constraint $G \succeq 0$ satisfies the complementarity condition $GM = 0$, so its columns immediately reveal such null directions.
In particular, for any global minimizer $x^*$, the Dirac measure $\delta_{x^*}$ is optimal for \cref{eq:gmp} and has a rank-one $\varphi$-moment matrix $M^* = \varphi(x^*)\varphi(x^*)^\top$ satisfying $G^*M^* = 0$.
Hence finite convergence enables exact recovery of each global minimizer from the preimage $\varphi^{-1}(v)$ of an associated rank-one dual optimum $M = vv^\top$.

Finite convergence is not required, however, to justify SOS hierarchies for numerical global optimization.
Asymptotic convergence already brings the evaluation vectors $\varphi(x^*)$ arbitrarily close (in Gramian seminorm) to satisfying the nullspace condition.
Thus, along any hierarchy, primal-dual solutions to \cref{eq:vsos} provide both increasingly tight lower bounds $f^*_\varphi$ and minimizer candidates of the form
\begin{equation}
    \hat{x} = \operatorname{proj}_\mathcal{V}(\tilde{x})\,,
    \qquad
    G\,\varphi(\tilde{x}) \approx 0
    \label{eq:extraction}
\end{equation}
Typically, $\tilde{x} \in \varphi^{-1}(v)$ with $M \approx vv^\top$, while $\operatorname{proj}_\mathcal{V}$ is a few steps of local optimization with discrete variables rounded.
Then the bound certifies the quality of $\hat{x}$ (or any feasible candidate) through the \textit{optimality gap} $\Delta_\varphi : \mathcal{V} \to \mathbb{R}_+$,
\begin{equation*}
    \Delta_\varphi(x) := f(x) - f^*_\varphi
    \;\;\geq\;\; f(x) - f^*
\end{equation*}
No feasible point can improve on $\hat{x}$ by more than $\Delta_\varphi(\hat{x})$, so once this gap falls within the application's tolerance, $\hat{x}$ is accepted as a numerical global optimum.

All that remains then is to select a particular hierarchy of polynomial bases $\varphi$ and to design a sufficient linear representation for membership in $\mathcal{I}(\mathcal{V})$.
With these choices, each instance of \cref{eq:vsos} becomes a concrete SDP that can be solved reliably by standard convex optimization algorithms.
In the next two subsections, we explain and compare the traditional choices with their sample-based counterparts.


\subsection{Traditional SOS}
\label{sec:traditional_sos}

The traditional choice of basis is all monomials up to some fixed degree $d \in \mathbb{N}$,
\begin{equation*}
    \varphi(x) \,=\, \nu_d(x) \,:=\, \big[x^\alpha\big]_{\alpha \in \mathbb{N}^n,\, |\alpha| \leq d}
\end{equation*}
Increasing $d$ yields a nested sequence of bases, and every polynomial eventually lies in its span.
However, its size grows combinatorially in problem dimension and degree,
\begin{equation*}
    N \,=\, \varpi(n,d) \,:=\, \binom{n+d}{d} \,=\, \frac{(n+d)!}{n!\,d!}
\end{equation*}

Next, to handle membership in $\mathcal{I}(\mathcal{V})$, the traditional approach replaces the full vanishing ideal with the ideal generated by the problem constraints,
\begin{equation*}
    \langle h_{1:\ell} \rangle \,:=\,
    \left\{
        \sum_{j=1}^{\ell} \lambda_j h_j
        \;\,\middle|\;\,
        \lambda_{1:\ell} \in \mathbb{R}[x_{1:n}]
    \right\}
\end{equation*}
Because each $h_j$ vanishes on $\mathcal{V}$ by construction, so does any polynomial combination of them.
Thus membership in $\langle h_{1:\ell} \rangle$ is sufficient for membership in $\mathcal{I}(\mathcal{V})$, though not necessary:
some polynomials that vanish on $\mathcal{V}$ may not be generated by $h_{1:\ell}$.\footnote{
    This follows from a foundational result in algebraic geometry known as the \textit{Real Nullstellensatz}:
    $\mathcal{I}(\mathcal{V})$ is the \textit{real radical} of $\langle h_{1:\ell} \rangle$, which can be strictly larger than $\langle h_{1:\ell} \rangle$ itself.
    As a simple example, if $\mathcal{V} = \{0\}$ is described by $h(x) = x^2$, then $p(x) = x$ vanishes on $\mathcal{V}$ but $p \neq \lambda h$ for any polynomial $\lambda$, i.e.\ $x \notin \langle x^2 \rangle$.
    All such missing polynomials are recovered by passing from $\langle x^2 \rangle$ to its radical $\langle x \rangle$.
}
\begin{equation*}
    \langle h_{1:\ell} \rangle \;\subseteq\; \mathcal{I}(\mathcal{V})
\end{equation*}
Lasserre showed in \cite{lasserre2001hierarchy} that replacing $\mathcal{I}(\mathcal{V})$ with $\langle h_{1:\ell} \rangle$ does not break convergence of the hierarchy.
It does, however, tie the rate of convergence to the particular algebraic description $h_{1:\ell}$ rather than the geometry of $\mathcal{V}$ itself.
Consequently, traditional SOS requires careful choice among seemingly equivalent problem formulations so that $\langle h_{1:\ell} \rangle$ captures more of $\mathcal{I}(\mathcal{V})$ at low degree.\footnote{
    Practitioners typically augment their problems with so-called \textit{redundant constraints}:
    additional polynomial equalities that do not change the feasible set but are included just to enlarge $\langle h_{1:\ell} \rangle$.
    The need to find such constraints is often treated as an inherent challenge in using SOS \cite{dumbgen2024autotight}, but it is really an artifact of approximating $\mathcal{I}(\mathcal{V})$ by $\langle h_{1:\ell} \rangle$ in the traditional approach.
}

With the monomial basis and generated ideal, the SDP \cref{eq:vsos} can be made concrete.
First, express the polynomial multipliers $\lambda_{1:\ell}$ on the same basis,
\begin{equation*}
    \lambda_j(x) \,=\, \varphi(x)^\top H_j\,\varphi(x)\,,
    \quad
    H_j \in \mathbb{S}^N
    \quad
    \forall j \in \{1, \ldots, \ell\}
\end{equation*}
where here $\varphi = \nu_d$ and $N = \varpi(n,d)$.\footnote{
    In most treatments, each multiplier uses a smaller monomial basis such that $\deg(\lambda_j h_j) \leq 2d$.
    We use the common-basis form here because it keeps the notation uniform and makes the connection to our approach in \cref{sec:circumventing_enumeration} clearer.
    This choice just admits some higher-degree multiplier products, which cannot weaken the relaxation for any fixed $d$.
}
Then the ideal membership condition becomes
\begin{equation*}
    \exists\,H_{1:\ell} \in \mathbb{S}^{N} :\;\; f - \beta - \textcolor{gray}{\overbrace{\normalcolor \varphi^\top G\,\varphi}^{\sigma}} \;=\; \sum_{j=1}^{\ell} \big(\textcolor{gray}{\overbrace{\normalcolor \varphi^\top H_j\,\varphi}^{\lambda_j}}\big) h_j
\end{equation*}
Note that unlike the Gramian $G$, the \textit{multiplier matrices} $H_{1:\ell}$ are not constrained to be PSD, because the polynomials $\lambda_{1:\ell}$ need not be SOS.
Introducing $H_{1:\ell}$ as decision variables and plugging in the monomial basis, the SDP becomes
\begin{equation}
    \max_{\substack{\beta \,\in\, \mathbb{R} \\ G, H_{1:\ell} \,\in\, \mathbb{S}^{\varpi(n,d)}}}\, \beta
    \quad \text{s.t.} \quad
    f - \beta - \nu_d^\top \bigg(G + \sum_{j=1}^{\ell} h_j H_j\bigg) \nu_d \,=\, 0\,,
    \quad
    G \succeq 0
    \label{eq:trad_sdp}
\end{equation}
The equality constrains the polynomial on the left side to be zero globally.
Traditional SOS enforces this by collecting like-terms and setting each collected coefficient to zero.
Since these collected coefficients are linear in the decision variables ($\beta$, $G$, and $H_{1:\ell}$), the ideal membership condition has become a finite linear system.

Thus for a fixed degree $d$, \cref{eq:trad_sdp} can be solved by a standard convex optimizer to obtain a certified lower bound and extract a solution candidate via \cref{eq:extraction}.
The hierarchy is simply to increase $d$ until the resulting optimality gap falls within tolerance.
The lowest meaningful value of $d$ is limited by the coefficient matching step used to construct the linear system.
A basic guideline is $d \geq \left\lceil \deg(f)/2 \right\rceil$ to ensure that the certificate can even match the leading terms of $f$.\footnote{
    Notably, if the problem polynomials have degree at most two, then the hierarchy can begin with $\nu_1 = [1,x]$, and the resulting SDP is the famous Shor relaxation for a quadratically constrained quadratic program (QCQP).
    Because that basis has just $N = n+1$ elements, practitioners often seek QCQP reformulations of higher-degree POPs before applying traditional SOS.
    However, this is largely an artifact of coefficient matching, which ties $d$ to the ambient degree of the objective even when its intrinsic degree on $\mathcal{V}$ is much lower.
}

\begin{figure}[htbp]
    \centering
    \includegraphics[width=0.75\linewidth]{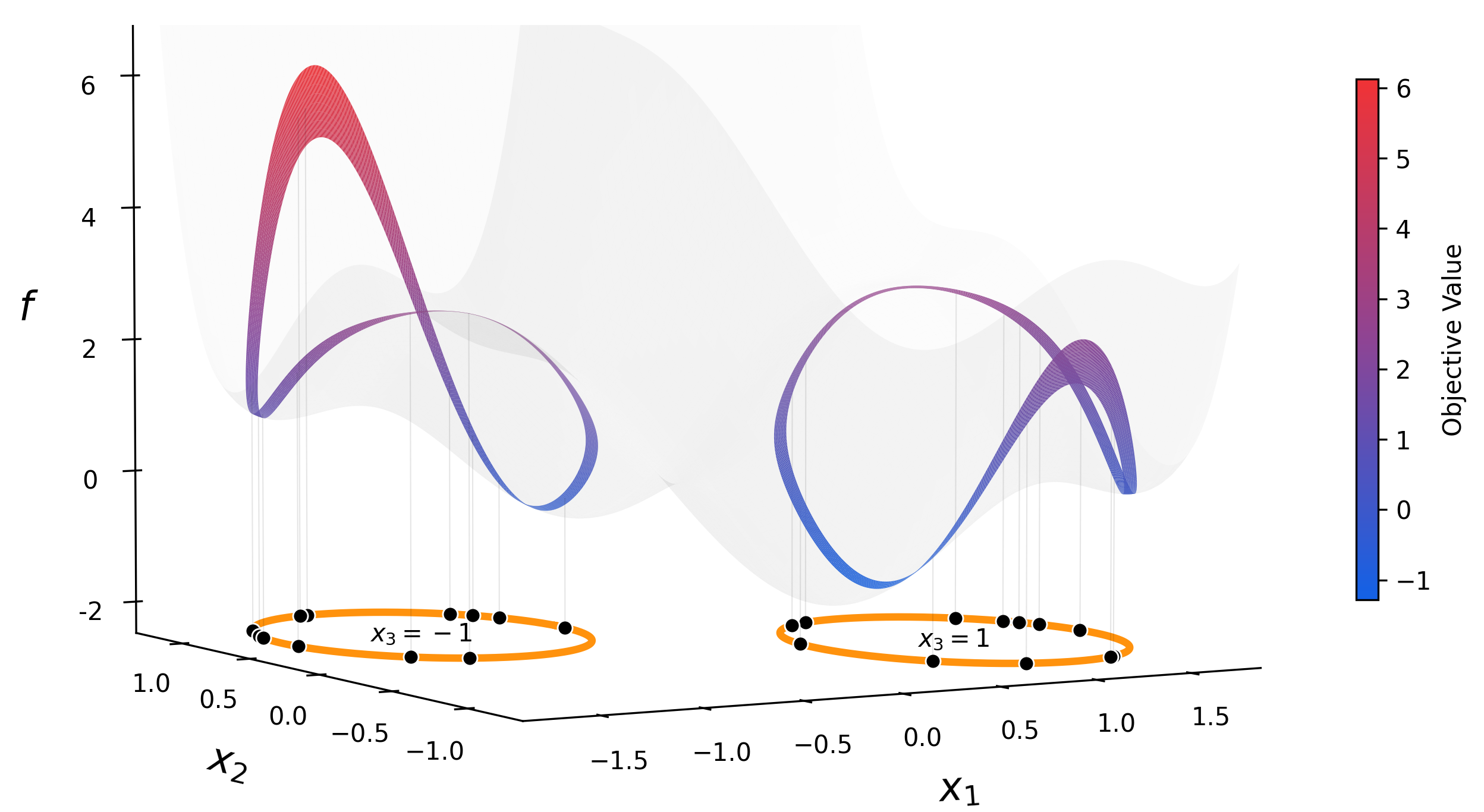}
    \caption{
    Visualization of the example POP \cref{eq:doublecircle}.
    The pair of orange circles is the feasible set projected onto the $(x_1,x_2)$-plane.
    The gray surface shows the objective, and the blue-red ribbons highlight its values on the feasible set.
    The samples used by sample-based SOS are marked in black.
    }
    \label{fig:doublecircle}
\end{figure}

To illustrate traditional SOS, we introduce the following simple POP that will also serve as a comparison point for sample-based SOS.
\begin{equation}
\begin{array}{r@{\quad}l}
    \displaystyle \min_{x \,\in\, \mathbb{R}^3}
    &
    \displaystyle \tfrac{1}{3}x_1^6 - 2x_1^4 + 4x_1^2 + x_1x_2 + 4x_2^4 + x_2^3 - 4x_2^2
    \\[8pt]
    \text{s.t.}
    &
    \vcenter{\hbox{$\begin{array}{r@{\;\;=\;\;}l}
        \left(x_1 - x_3\right)^2 + \left(x_2 + \tfrac{1}{2}x_3\right)^2 - \tfrac{1}{2} & 0 \\[4pt]
        x_3^2 - 1 & 0
    \end{array}$}}
\end{array}
\label{eq:doublecircle}
\end{equation}
The second constraint is equivalent to $x_3 \in \{-1,1\}$, so it switches the first constraint between two circles in $(x_1,x_2)$-space (at $(-1,1/2)$ and $(1,-1/2)$ respectively).
Thus the problem is to simultaneously decide which circle contains the minimizer and where the minimizer is on that circle.
The objective is visualized in \cref{fig:doublecircle}.

Since $\deg(f) = 6$, the first level of traditional SOS uses $d = \lceil 6/2 \rceil = 3$.
With $n = 3$ variables, the corresponding monomial basis has $N = \varpi(3,3) = 20$ elements.
One possible ordering is
\begin{equation*}
\begin{aligned}
    \nu_3(x) =
    \big[&
    1,\;
    x_1,\;x_2,\;x_3,\;
    x_1^2,\;x_1x_2,\;x_1x_3,\;x_2^2,\;x_2x_3,\;x_3^2, \\
    &
    x_1^3,\;x_1^2x_2,\;x_1^2x_3,\;
    x_1x_2^2,\;x_1x_2x_3,\;x_1x_3^2,\;
    x_2^3,\;x_2^2x_3,\;x_2x_3^2,\;x_3^3
    \big]^\top
\end{aligned}
\end{equation*}
Then to instantiate \cref{eq:trad_sdp}, traditional SOS symbolically expands the polynomial
\begin{equation*}
    f - \beta - \nu_3^\top \big(G + h_1H_1 + h_2H_2\big) \nu_3
\end{equation*}
and sets its coefficients to zero, resulting in one linear constraint per term. E.g.,
\begin{equation*}
\begin{aligned}
    \text{[1]}:\quad & -\beta - G_{1,1} + \tfrac{1}{2}(H_1)_{1,1} + (H_2)_{1,1} \,=\, 0 \\[2pt]
    \text{[$x_1x_2$]}:\quad & 1 - 2G_{1,6} - 2G_{2,3} + (H_1)_{1,6} + (H_1)_{2,3} + 2(H_2)_{1,6} + 2(H_2)_{2,3} \,=\, 0 \\[2pt]
    \text{[$x_1^6$]}:\quad & \tfrac{1}{3} - G_{11,11} - 2(H_1)_{2,11} - (H_1)_{5,5} + \tfrac{1}{2}(H_1)_{11,11} + (H_2)_{11,11} \,=\, 0
\end{aligned}
\end{equation*}
After eliminating redundant multiplier coefficients (here, any that produce terms above degree $6$), the resulting linear system has $\varpi(3,6) = 84$ constraints in the variables $\beta$, $G$, and the $70$ remaining scalar degrees of freedom in $H_{1:2}$.

Thus for this example, \cref{eq:trad_sdp} has one $20 \times 20$ PSD matrix variable and $71$ additional scalar variables.
Solving this SDP with the commercial convex optimizer \href{https://www.mosek.com/}{MOSEK} yields the lower bound $\beta = -1.239001081089$ and a Gramian $G$ with numerical corank $1$.
Taking its null vector $v$ scaled so that $v_1 = 1$, we recover \mbox{$\hat{x} = \tilde{x} = \nu_3^{-1}(v) = (v_2,v_3,v_4) \approx (0.352745, -0.784712, 1)$}.
This numerically feasible point has cost within $10^{-13}$ of the bound, indicating global optimality.

Such exactness on the first level of the hierarchy is common in practice, but even the first level can be intractable at realistic problem scales.
With the same modest $d = 3$, increasing to just $n = 100$ variables would require a Gramian of size $\varpi(100,3) = 176851$, already far beyond the capabilities of modern SDP solvers.
True scalability therefore requires exploiting additional problem structure such as sparsity.
Many well-studied extensions of traditional SOS do exactly that, but they are still built from overparameterized monomial bases and inherit its sensitivity to an algebraic description of $\mathcal{V}$.
We instead turn to sample-based SOS, which constructs the hierarchy from polynomial evaluations on $\mathcal{V}$ rather than hand-crafted equations.

\subsection{Sample-Based SOS}
\label{sec:sample_based_sos}

Global nonnegativity of the certificate $\sigma = \varphi^\top G\,\varphi$ is guaranteed by $G \succeq 0$ for any choice of basis $\varphi$.
The basis just determines the freedom available for $\sigma$ to match $f-\beta$.
Since this match is only required on $\mathcal{V}$, basis functions that are equivalent on $\mathcal{V}$ serve the same purpose and can be merged without weakening the hierarchy.
For example, $\mathcal{V}$ in \cref{eq:doublecircle} satisfies $x_3^2 = 1$ and thus the basis function $x_3^2$ in $\nu_3$ offers no expressivity beyond what $1$ already provides.
Many more such redundancies are far less obvious.
Samples from $\mathcal{V}$ furnish a numerical method to efficiently discover all of them.

\looseness=-1
Let $\mathcal{Z} := \{z_1, z_2, \ldots, z_S\} \subset \mathcal{V}$ be a finite set of feasible samples.
Evaluate all monomials up to degree $d$ at each sample to form the \textit{monomial evaluation matrix},
\begin{equation*}
    E_d(\mathcal{Z}) \,:=\,
    \begin{bmatrix}
        \nu_d(z_1)^\top \\
        \vdots \\
        \nu_d(z_S)^\top
    \end{bmatrix}
    \,\in\, \mathbb{R}^{S \times \varpi(n,d)}
\end{equation*}
Every null vector $c \in \ker E_d(\mathcal{Z})$ defines a polynomial $p = c^\top \nu_d$ that vanishes on all the samples.
Thus each degree-$\leq\!d$ relation that holds on $\mathcal{Z}$ corresponds to a linear dependence among the columns of $E_d(\mathcal{Z})$.
Since $\mathcal{Z} \subset \mathcal{V}$, any relation valid on all of $\mathcal{V}$ must also appear among the dependencies observed on $\mathcal{Z}$.
The only concern is whether the samples introduce additional relations that do not hold on $\mathcal{V}$.
Fortunately, finite-dimensional function spaces enable finite sampling to rule out such spurious relations.

\begin{definition}
    \label{def:d_poised}
    Sample set $\mathcal{Z}$ is \textit{$d$-poised for} variety $\mathcal{V}$ if, for any $p \in \operatorname{span}(\nu_d)$, vanishing on the samples requires vanishing on the whole variety:
    \begin{equation*}
        p(x) = 0 \;\ \forall x \in \mathcal{Z}
        \quad \implies \quad
        p(x) = 0 \;\ \forall x \in \mathcal{V}
    \end{equation*}
\end{definition}
I.e., $c^\top\!\nu_d \in \mathcal{I}(\mathcal{Z}) \,\Rightarrow\, c^\top\!\nu_d \in \mathcal{I}(\mathcal{V})$.
Geometrically, this means that the samples do not accidentally all lie on some subvariety of $\mathcal{V}$ cut out by extra degree-$\leq\!d$ relations.

Assuming $\mathcal{Z}$ is $d$-poised for $\mathcal{V}$, we can use the singular value decomposition (SVD) of $E_d(\mathcal{Z})$ to compute a basis $\varphi$ with no degree-$\leq\!d$ redundancies.
Let $r := \operatorname{rank} E_d(\mathcal{Z})$ and write its thin SVD as $E_d(\mathcal{Z}) = U_r \Sigma_r V_r^\top$.
The sample-based \textit{quotient basis}\footnote{
    Formally, this basis spans the degree-$\leq\!d$ portion of the \textit{coordinate ring} $\mathbb{R}[x_{1:n}] \,\big/\, \mathcal{I}(\mathcal{V})$.
    The quotient identifies polynomials that behave the same on $\mathcal{V}$, with the basis functions as representatives.
}
is
\begin{equation*}
    \varphi(x) \,=\, B^\top \nu_d(x)
    \,,\qquad
    B \,:=\, V_r \Sigma_r^{-1} \,\in\, \mathbb{R}^{\varpi(n,d) \times r}
\end{equation*}
Its evaluation matrix $E_d(\mathcal{Z})B = U_r$ is orthonormal (which also improves numerical conditioning).
Viewed as principal component analysis of the monomial ``features'' $\nu_d$ with uncentered covariance $E_d(\mathcal{Z})^\top E_d(\mathcal{Z})$, the SVD discards directions with no useful variation over the ``data'' $\mathcal{Z}$ and orthonormalizes the rest.
Poisedness guarantees that $\varphi$ preserves the expressivity of $\nu_d$ on $\mathcal{V}$ using only $N = r \leq \varpi(n,d)$ basis functions.

The same samples can also be used to enforce the ideal membership in \cref{eq:vsos}.
Let
\begin{equation*}
    \deg_{\mathcal{V}}(f) \,:=\, \min \big\{\deg(p)\ \big|\ f - p \in \mathcal{I}(\mathcal{V})\,\big\}
    \,\leq\, \deg(f)
\end{equation*}
be the \textit{effective degree} of $f$ on $\mathcal{V}$.
If $\deg_{\mathcal{V}}(f) \leq 2d$ and $\mathcal{Z}$ is $2d$-poised for $\mathcal{V}$, then $f$ has a degree-$\leq\!2d$ representative on $\mathcal{V}$ and poisedness immediately gives,
\begin{equation*}
    f - \beta - \varphi^\top G\,\varphi \;\in\; \mathcal{I}(\mathcal{V})
    \quad \iff \quad
    f(x) - \beta - \varphi(x)^\top G\,\varphi(x) \,=\, 0 \ \ \forall x \in \mathcal{Z}
\end{equation*}
Simply put, we can match $\sigma$ to $f-\beta$ on $\mathcal{V}$ by enforcing that they match at the samples.
This produces one linear equation per sample and defines the sample-based~SDP,
\begin{equation}
    \max_{\substack{\beta \,\in\, \mathbb{R} \\ G \,\in\, \mathbb{S}^r_+}}\, \beta
    \quad \text{s.t.} \quad
    f(x) - \beta \,=\, \varphi(x)^\top G\,\varphi(x) \,\ \ \forall x \in \mathcal{Z}
    \label{eq:samp_sdp}
\end{equation}
Letting $f(\mathcal{Z}) := [f(z_1), \ldots, f(z_S)]^\top$ and defining the \textit{product evaluation matrix}
\begin{equation*}
    \Phi(\mathcal{Z}) \,:=\,
    \begin{bmatrix}
        \operatorname{svec}\!\big(\varphi(z_1)\varphi(z_1)^\top\big)^\top \\
        \vdots \\
        \operatorname{svec}\!\big(\varphi(z_S)\varphi(z_S)^\top\big)^\top
    \end{bmatrix}
    \,\in\, \mathbb{R}^{S \times \tfrac{r(r+1)}{2}}
\end{equation*}
we can neatly arrange the linear system in \cref{eq:samp_sdp} as
\begin{equation*}
    \begin{bmatrix}
        \mathbf{1}_S & \Phi(\mathcal{Z})
    \end{bmatrix}
    \begin{bmatrix}
        \beta \\
        \operatorname{svec}(G)
    \end{bmatrix}
    \,=\, f(\mathcal{Z})
\end{equation*}

We can apply sample-based SOS to \cref{eq:doublecircle} using just $S=24$ samples,\footnote{
    Here $24$ is the dimensionality of the space of degree-$\leq\!6$ polynomials on $\mathcal{V}$.
    \Cref{sec:selecting_samples} shows how this required sample count can be determined numerically.
}
evenly split between the two circles as shown in \cref{fig:doublecircle}.
The SVD of $E_3(\mathcal{Z})$ reveals rank $r=12$, reducing the traditional $20$-element monomial basis to $12$ better-conditioned basis polynomials.
Hence \cref{eq:samp_sdp} has one $12 \times 12$ PSD matrix variable $G$, one scalar variable $\beta$, and no multipliers.
Whereas traditional SOS produced an $84 \times 281$ linear system through symbolic coefficient matching, the sampled system $\big([\mathbf{1}_{24}\ \Phi(\mathcal{Z})],\,f(\mathcal{Z})\big)$ is only $24 \times 79$ and was built entirely from black-box evaluations of $\nu_3$ and $f$.

Furthermore, poisedness makes sample-based SOS at least as tight as traditional SOS because $\mathcal{I}(\mathcal{V})$ contains every possible description $\langle h_{1:\ell} \rangle$ that traditional SOS could use.
Thus for \cref{eq:doublecircle}, solving \cref{eq:samp_sdp} recovers the same tight bound as before.
Solution extraction can likewise proceed from a null vector $v$ of $G$:
poisedness gives $\nu_d = V_r \Sigma_r \varphi$ on $\mathcal{V}$, so normalizing $V_r \Sigma_r v$ recovers the monomial coordinates that expose $\tilde{x}$.

Sample-based SOS also offers a basis-agnostic route to extraction that will be very convenient later.
Consider the dual of \cref{eq:samp_sdp} with one variable per sample equality,
\begin{equation*}
    \min_{y \,\in\, \mathbb{R}^S}\,\, f(\mathcal{Z})^\top y
    \quad\,\ \text{s.t.}\ \ \quad
    \mathbf{1}_S^\top y = 1
    \,,\,\quad
    \sum_{s=1}^{S} y_s \varphi(z_s)\varphi(z_s)^\top \succeq 0
\end{equation*}
The \textit{sample weights} $y_{1:S}$ define the normalized linear functional
\begin{equation*}
    \mathcal{L}_y[p] \,:=\, \sum_{s=1}^{S} y_s p(z_s)
    \,,\ \ \quad
    \mathcal{L}_y[1] = 1
\end{equation*}
whose moments play the same role as those of the measure $\mu$ in \cref{eq:gmp}.\footnote{
    Note that $y_{1:S}$ are signed and cannot themselves be interpreted as probabilities.
    It is their induced moments that correspond to those of $\mu$ and support solution extraction.
}
This lets us query any desired dual information directly, regardless of the basis used to construct the SDP.
In particular, the first-moment vector
\begin{equation*}
    \bar{x}
    \,:=\,
    \mathcal{L}_y[x]
    \,=\,
    \sum_{s=1}^{S} y_s z_s
\end{equation*}
coincides with $\tilde{x}$ when $G$ has corank $1$ and otherwise provides a natural solution candidate when no null direction is distinguished.\footnote{
    Higher-order moments can further aid extraction, for example by revealing the spread of multiple global optima or by conditioning the mean on likely discrete branches.
}
Indeed for \cref{eq:doublecircle}, this affine combination of the samples pinpoints $x^*$.

Overall, sample-based SOS replaces the traditional relaxation with a stronger yet smaller SDP that is invariant to the chosen constraints $h_{1:\ell}$ and constructed entirely from black-box evaluations.
Despite these advantages, the formulation as posed still faces several obstacles to practical scale.
First, it assumes access to a sufficiently poised sample set, which we will see can become intractable for problems with many discrete variables.
Second, the current equalities-only ($\mathcal{F} = \mathcal{V}$) treatment handles inequalities via slack variables that bloat the ambient dimension.
Third, and perhaps most importantly, this formulation does not naturally expose sparsity: both $E_d(\mathcal{Z})$ and $\Phi(\mathcal{Z})$ are large dense matrices even when the original problem is sparse.
Fourth, because each sample becomes an SDP constraint, we would like a way to detect sufficiency and discard the excess before solving.

Our four contributions address these issues in turn.
We then demonstrate their combined effect on two real-world engineering problems.
This methodology reaches unprecedented scales for certifiable semialgebraic optimization using only black-box evaluations and an off-the-shelf SDP solver.


\section{Circumventing Enumeration}
\label{sec:circumventing_enumeration}

Randomization provides a natural way to achieve poisedness.
By continuity, every polynomial equality that fails somewhere on $\mathcal{V}$ also fails on an open subset of $\mathcal{V}$.
Thus under any distribution fully supported on $\mathcal{V}$, a random sample has positive probability of landing there and ruling out the equality.
Therefore, independent random draws progressively eliminate spurious degree-$\leq\!d$ relations, and the probability that $\mathcal{Z}$ is $d$-poised approaches one.
Moreover, random feasible sampling is already commonplace in applications, e.g.\ simulation rollouts for trajectory optimization or evolutionary algorithms for design optimization.

The difficulty arises when $\mathcal{V}$ is \textit{reducible}, i.e.\ when it is a finite union of proper subvarieties.
Then the sampling distribution must spread its probability mass across these components, which can leave some unlikely to be sampled.
Missing even a single component can cause the samples to satisfy relations that are false on $\mathcal{V}$ as a whole.
For instance, if the samples in \cref{fig:doublecircle} landed only on the $x_3=-1$ circle, then $x_3+1=0$ would appear to hold on $\mathcal{V}$ and \cref{eq:samp_sdp} would become blind to the other circle and true solution.
The componentwise strategy of \cite{cifuentes2017sampling} addresses this issue by enumerating the irreducible components of $\mathcal{V}$ and sampling each separately.
However, discrete variables can fragment $\mathcal{V}$ into exponentially many components, making enumeration intractable even for modest-sized mixed-integer problems.

We propose instead to sample a larger variety $\mathcal{W} \supseteq \mathcal{V}$ whose geometry is easier to cover, while retaining any violated equality constraints through traditional multipliers.
To see this, let $\mathcal{Z} \subset \mathcal{W}$ and construct the quotient basis $\varphi = B^\top \nu_d$ from $E_d(\mathcal{Z})$ just as in \cref{sec:sample_based_sos}.
Then augment \cref{eq:samp_sdp} with one polynomial multiplier $\lambda_j = \varphi^\top H_j\,\varphi$ for each original equality to form the \textit{ambient-variety SDP},
\begin{equation}
    \beta^*_\mathcal{W} \,:=
    \max_{\substack{\beta \,\in\, \mathbb{R} \\ G \,\in\, \mathbb{S}^r_+ \\ H_{1:\ell} \,\in\, \mathbb{S}^r}}\, \beta
    \quad \text{s.t.} \quad
    f(x) - \beta \,=\, \varphi(x)^\top \bigg(G + \sum_{j=1}^\ell h_j(x)H_j\bigg)\, \varphi(x)
    \;\ \forall x \in \mathcal{Z} \subset \mathcal{W}
\label{eq:ambi_sdp}
\end{equation}
Since every $h_j$ vanishes on $\mathcal{V}$, choosing $\mathcal{W}=\mathcal{V}$ eliminates all multiplier terms and reduces \cref{eq:ambi_sdp} to \cref{eq:samp_sdp}.
At the other extreme, choosing $\mathcal{W}=\mathbb{R}^n$ makes \cref{eq:ambi_sdp} equivalent to \cref{eq:trad_sdp} under the poisedness condition formalized below.
Specifically, the basis would have $r=\varpi(n,d)$ and matching at the samples would just replace coefficient matching as an equivalent way to enforce the same polynomial identity.
Even then, orthonormality and black-box assembly are preserved as in \cite{lofberg2004sampling}.
Choosing $\mathcal{W}$ between these extremes essentially interpolates between \cref{eq:samp_sdp} and \cref{eq:trad_sdp}, with the resulting bounds ordered by the following theorem.

\begin{theorem}[Ambient-Variety Bounds]
    \label{thm:ambi_bounds}
    Fix $d$ and let $\beta^*_\mathcal{V}$ and $\beta^*_{\mathbb{R}^n}$ denote the optimal values of \cref{eq:samp_sdp} and \cref{eq:trad_sdp} respectively.
    Suppose $\mathcal{Z}$ is $\kappa$-poised for $\mathcal{W} \supseteq \mathcal{V}$, with $\kappa \geq \max\{2d,\deg_\mathcal{W}(f)\}$ and $\kappa \geq 2d+\deg_\mathcal{W}(h_j)$ for each constraint polynomial $h_j$ that does not vanish on $\mathcal{W}$.
    Then,
    \begin{equation*}
        \beta^*_{\mathbb{R}^n}
        \;\leq\;
        \beta^*_\mathcal{W}
        \;\leq\;
        \beta^*_\mathcal{V}
        \;\leq\;
        f^*
    \end{equation*}
\end{theorem}
\begin{proof}
    Let $p := f - \beta - \sigma - \sum h_j\lambda_j$.
    Under the stated poisedness condition, the sampled system in \cref{eq:ambi_sdp} is equivalent to $p \in \mathcal{I}(\mathcal{W})$.
    Since $\mathcal{V} \subseteq \mathcal{W} \subseteq \mathbb{R}^n$, their vanishing ideals satisfy $\mathcal{I}(\mathbb{R}^n) = \{0\} \subseteq \mathcal{I}(\mathcal{W}) \subseteq \mathcal{I}(\mathcal{V})$.
    Thus every certificate feasible for \cref{eq:trad_sdp} can be reparameterized on the quotient basis and is feasible for \cref{eq:ambi_sdp}, giving $\beta^*_{\mathbb{R}^n} \leq \beta^*_\mathcal{W}$.
    Since both $p$ and $\sum h_j\lambda_j$ belong to $\mathcal{I}(\mathcal{V})$, so does their sum $f - \beta - \sigma$.
    Thus any $\sigma$ in \cref{eq:ambi_sdp} can certify the same bound in \cref{eq:samp_sdp}, giving $\beta^*_\mathcal{W} \leq \beta^*_\mathcal{V}$.
\end{proof}

The same ideal-inclusion argument shows that shrinking $\mathcal{W}$ through nested varieties can only strengthen the bound.
Moreover, \cref{thm:ambi_bounds} transfers the convergence of traditional SOS on compact $\mathcal{V}$ to the entire ambient-variety hierarchy obtained by increasing $d$, provided that the poisedness condition is maintained at each level.
Importantly, the chosen ambient variety $\mathcal{W}$ itself need not be compact.

When poised sampling of $\mathcal{V}$ is difficult, it's typically because a small subset of the constraints $h_{1:\ell}$ fragment $\mathcal{V}$ into a prohibitively large number of components, while the other constraints just define geometry that is simple to cover.
We therefore propose to build $\mathcal{W} \supseteq \mathcal{V}$ using only those tractable constraints.
Let $J \subseteq \{1,\ldots,\ell\}$ index the equalities to be retained through explicit multipliers, and define the ambient variety using the remaining equalities,
\begin{equation*}
    \mathcal{W} \,:=\, \Big\{x \in \mathbb{R}^n \;\,\Big\vert\;\, h_j(x)=0 \;\ \forall j \in \{1,\ldots,\ell\} \setminus J\Big\}
\end{equation*}
The multiplier terms with $j \notin J$ vanish on $\mathcal{W}$ and drop out of \cref{eq:ambi_sdp}.
Thus only the equalities indexed by $J$ incur multiplier matrices, while those respected by $\mathcal{W}$ serve to reduce the polynomial basis.\footnote{
    This quotient basis spans the degree-bounded portion of the coordinate ring $\mathbb{R}[x_{1:n}] \,\big/\, \mathcal{I}(\mathcal{W})$.
    The respected equalities enlarge $\mathcal{I}(\mathcal{W})$ and thus quotient out more polynomial directions.
}

To make this idea concrete, return to example \cref{eq:doublecircle} and let $J=\{2\}$.
Omitting the binary constraint $h_2$ from $\mathcal{W}$ merges the two distinct circles into one irreducible surface $h_1(x)=0$.
Then, instead of enumerating $x_3=\pm1$, we can draw $x_3$ from a normal distribution and sample $(x_1,x_2)$ on the resulting circular slice.
The multiplier for $h_2$ upholds the binary constraint in the certification, so the resulting $d=3$ relaxation recovers the same tight bound and solution as before.
However, the basis reduction is less extensive, yielding $r=16$ functions rather than the $r=12$ obtained from samples poised for $\mathcal{V}$ itself.
Thus we have exchanged some quotient efficiency for much easier sampling.
For problems with many discrete variables, this approach avoids combinatorial branch enumeration while preserving basis reductions from the continuous geometry and a bound at least as tight as traditional SOS.


\section{Handling Inequalities}
\label{sec:handling_inequalities}

So far, we have worked in the equalities-only setting $\mathcal{F}=\mathcal{V}$ by treating the $m$ inequalities $g_i(x) \geq 0$ as equalities $g_i(x) - s_i^2 = 0$ with slacks $s_{1:m}$ absorbed into $x$.
Such lifting is exact, but it inflates the degree-$d$ monomial count from $\varpi(n,d)$ to $\varpi(n+m,d)$ and can result in a much larger SDP even after quotienting.
We now return to the original feasible set $\mathcal{F} = \{x \in \mathcal{V} \mid g_{1:m}(x) \geq 0\}$ and show how to handle the inequalities directly.

Since each $g_i$ is nonnegative on $\mathcal{F}$, its product with any SOS polynomial $\sigma_i$ is also nonnegative there.
Thus for any SOS polynomials $\sigma_{0:m}$, the combination $\sigma_0 + \sum_{i=1}^{m} g_i\sigma_i$ is nonnegative on $\mathcal{F}$ and can serve as a certificate by matching $f-\beta$ on $\mathcal{V} \supseteq \mathcal{F}$.
Construct a basis $\varphi = B^\top \nu_d$ from samples $\mathcal{Z} \subset \mathcal{W}$ as in \cref{sec:circumventing_enumeration}, but augment the $\kappa$-poisedness condition of \cref{thm:ambi_bounds} with $\kappa \geq 2d+\max_i\deg_\mathcal{W}(g_i)$.
Then writing $\sigma_i = \varphi^\top G_i\,\varphi$ and $\lambda_j = \varphi^\top H_j\,\varphi$ gives the SDP,
\begin{equation}
    \max_{\substack{\beta \,\in\, \mathbb{R} \\ G_{0:m} \,\in\, \mathbb{S}^r_+ \\ H_J \,\in\, \mathbb{S}^r}}\, \beta
    \quad \text{s.t.} \quad
    f(x) - \beta \,=\, \varphi(x)^\top \bigg(G_0 + \sum_{i=1}^{m} g_i(x)G_i + \sum_{j \in J} h_j(x)H_j\bigg)\, \varphi(x)
    \;\ \forall x \in \mathcal{Z}
    \label{eq:ineq_sdp}
\end{equation}
Note that the samples themselves need not satisfy the inequalities.
And as before, each sample yields one linear equation via black-box evaluation of the model polynomials.

To establish convergence of the optimal bounds as $d$ is increased, we use the following theorem adapted from \cite[Section~3.6.2]{laurent2008book}.
\begin{theorem}[Putinar's Positivstellensatz]
    \label{thm:putinar}
    Suppose that the boundedness of $\mathcal{F}$ is SOS-certified by the constraints defining it.
    If a real polynomial $p$ is strictly positive on $\mathcal{F}$, then there exist SOS polynomials $\sigma_{0:m}$ such that $p - \big(\sigma_0 + \sum_{i=1}^{m} g_i\sigma_i\big) \,\in\, \langle h_{1:\ell} \rangle$.
\end{theorem}
The first hypothesis is the \textit{Archimedean property}, a slightly stronger version of our standing compactness assumption that is easily ensured in practice by including a norm bound $b-\|x\|^2\geq0$ among $g_{1:m}$.
Then for any $\beta < f^*$, \cref{thm:putinar} guarantees a representation of $f-\beta$ on $\mathcal{W}$ that is feasible for \cref{eq:ineq_sdp} at sufficiently large $d$.
Since these certificate families are nested, their optimal bounds converge monotonically to $f^*$ as $d$ increases, provided that the stated poisedness conditions hold at each level.
The dual sample weights still provide solution candidates as in \cref{sec:sample_based_sos}.

Moreover, this formulation exposes how the bound can be tightened at fixed $d$ by incorporating products of distinct inequality polynomials.
Since any product of $g_{1:m}$ is also nonnegative on $\mathcal{F}$, each can contribute a certificate term such as $g_i g_j g_k \sigma_{ijk}$ with its own SOS multiplier.
Using all such terms up to a fixed degree (with poisedness adjusted accordingly) gives a sample-based analogue of the \textit{Schm\"udgen hierarchy}.

The SDP \cref{eq:ineq_sdp} has $m+1$ PSD blocks of size $r$.\footnote{
    Uniform block size follows from using the same basis $\varphi$ for every SOS polynomial.
    Another natural choice is to impose $\deg(g_i \sigma_i) \leq 2d$, allowing smaller $G_i$ at the risk of a weaker bound.
}
To compare this with slack lifting, consider the inequalities-only setting $\ell=0$ and assume that the interior of $\mathcal{F}$ is nonempty.
Even after quotienting by all slack equalities $s_i^2=g_i(x)$, the slack-lifted SDP has a single PSD block of size
\begin{equation}
    r_{\text{slack-lifted}} \,\geq\,
    \sum_{k=0}^{\min(m,d)} \binom{m}{k} \varpi(n,d-k)
    \,=\, \frac{m^d}{d!} + \mathcal{O}\!\left(m^{d-1}\right)
    \label{eq:slack_bound}
\end{equation}
The sum counts monomials of the form $s_{i_1}\!\cdots s_{i_k} x^\alpha$ with distinct slack indices and $|\alpha| \leq d-k$.
Those monomials remain independent on the lifted variety because terms containing different subsets of slacks cannot cancel for every choice of signs in $s_i=\pm\sqrt{g_i(x)}$, and the $x^\alpha$ are independent on the feasible interior.
Thus for fixed $n$ and $d$, \cref{eq:slack_bound} shows that the slack-lifted PSD block size grows no slower than $m^d$, whereas our \cref{eq:ineq_sdp} grows linearly in $m$ across fixed-size PSD blocks.


\section{Leveraging Sparsity}
\label{sec:leveraging_sparsity}

Most real large-scale optimization problems exhibit significant structural sparsity.
For polynomial models, this sparsity can manifest in two distinct forms.
\textit{Term sparsity} occurs when each polynomial contains only a small fraction of the monomials available at its degree, allowing a smaller SOS basis to be determined through combinatorial analysis of the exponent geometry.
However, exploiting term sparsity first requires identifying those monomials, which would undermine our goal of operating directly through black-box evaluations.
We instead focus on \textit{correlative sparsity}, which occurs when each polynomial depends on only a small subset of the decision variables.
This is the same locality captured by factor graphs in e.g.\ planning and sensor fusion, where each factor couples only a few variables even when the full problem is very large.

Suppose that the objective function is given as
\begin{equation*}
    f \,=\, \sum_{k=1}^{K} f_k
    \,,\qquad
    f_k \,\in\, \mathbb{R}[x_{C_k}]
\end{equation*}
where the \textit{cliques} $C_1,\ldots,C_K \subseteq \{1,\ldots,n\}$ index the decision variables $x_{C_k}$ on which each \textit{objective component} $f_k$ may depend.
Likewise, let $\mathcal{G}_{1:K}$ partition $\{1,\ldots,m\}$ and $\mathcal{H}_{1:K}$ partition $\{1,\ldots,\ell\}$ such that
\begin{equation*}
    g_i \,\in\, \mathbb{R}[x_{C_k}] \quad \forall i \in \mathcal{G}_k
    \,,\qquad
    h_j \,\in\, \mathbb{R}[x_{C_k}] \quad \forall j \in \mathcal{H}_k
\end{equation*}
In other words, each objective component and constraint is assigned to a clique that contains its \textit{scope}, the decision variables on which it explicitly depends.
The cliques may overlap, so a decision variable can participate in several parts of the model.

To make this concrete, consider the following toy problem:
\begin{equation}
\begin{aligned}
    \min_{x \,\in\, \mathbb{R}^5}\quad
    & -x_2x_3 + (x_4^2 - x_5)
    \\[4pt]
    \text{s.t.}\quad
    &
    \begin{alignedat}{2}
        g_1(x_2) &\,:=\, x_2 - 2x_2^2 \,\geq\, 0
        &\quad\ g_2\smash[t]{\textcolor{gray}{\overbrace{\normalcolor (x_3,x_4,x_5)}^{\text{scope}}}} &\,:=\, 1 - x_3^2 - x_4^2 - x_5^2 \,\geq\, 0
        \\[3pt]
        h_1(x_1,x_2) &\,:=\, x_1 - x_2^2 \,=\, 0
        &\quad\ h_2(x_2,x_3) &\,:=\, x_2^2 + x_3^2 - 1 \,=\, 0
    \end{alignedat}
\end{aligned}
\label{eq:csp_example}
\end{equation}
The scopes of the objective components and constraints determine the problem's \textit{correlative sparsity pattern} (CSP), represented by a graph where two decision variables are connected by an edge if they occur together in any scope.
Here, the maximal cliques of this graph provide a natural decomposition:
\begin{equation*}
    \vcenter{\hbox{
        \begin{tikzpicture}[
            x=8mm,
            y=8mm,
            line width=0.4pt,
            every node/.style={inner sep=1.5pt}
        ]
            \node (v1) at (0,1)   {$x_1$};
            \node (v2) at (1.3,1) {$x_2$};
            \node (v3) at (2.6,1) {$x_3$};
            \node (v4) at (1.9,0) {$x_4$};
            \node (v5) at (3.3,0) {$x_5$};
            \draw (v1) -- (v2) -- (v3);
            \draw (v3) -- (v4) -- (v5) -- (v3);
        \end{tikzpicture}
    }}
    \qquad
    \begin{array}{c|c|c|c|c}
        k &       C_k &       f_k & \mathcal{G}_k & \mathcal{H}_k \\
        \hline
        1 & \{1,2\}   & 0         &         \{1\} &         \{1\} \\
        2 & \{2,3\}   & -x_2x_3   &         \{\}  &         \{2\} \\
        3 & \{3,4,5\} & x_4^2-x_5 &         \{2\} &         \{\}
    \end{array}
\end{equation*}

One can always take the trivial decomposition with $K=1$ and $C_1=\{1,\ldots,n\}$, but correlative sparsity becomes useful when the problem admits many small cliques with even smaller overlaps, analogous to low treewidth in graph theory.
Of particular importance is the \textit{running intersection property} (RIP), which holds if the cliques can be ordered such that, for every $k>1$, there exists a parent index $\pi(k)<k$ satisfying
\begin{equation*}
    C_k \cap \bigcup_{l<k} C_l \,\subseteq\, C_{\pi(k)}
    \tag{RIP}
\end{equation*}
I.e., the cliques can be arranged as a tree such that all cliques sharing any particular variable form a connected subtree.
Our decomposition of \cref{eq:csp_example} satisfies RIP in the displayed order, with $\pi(2)=1$ and $\pi(3)=2$.
The RIP plays a key role in the following sparse variant of Putinar's Positivstellensatz, adapted from \cite[Theorem~5]{grimm2007sparsity}.

\begin{theorem}[Correlatively Sparse Positivstellensatz]
    \label{thm:spositiv}
    Let $C_{1:K}$ satisfy RIP and suppose that the constraints assigned to each $C_k$ SOS-certify boundedness of $x_{C_k}$ using polynomials in $\mathbb{R}[x_{C_k}]$.
    If a real polynomial $p = \sum_{k=1}^{K} p_k$ with $p_k \in \mathbb{R}[x_{C_k}]$ is strictly positive on $\mathcal{F}$, then there exist SOS polynomials $\sigma_{k,i}$ for $i \in \{0\} \cup \mathcal{G}_k$ and unrestricted polynomials $\lambda_{k,j}$ for $j \in \mathcal{H}_k$, all in $\mathbb{R}[x_{C_k}]$, such that
    \begin{equation*}
        p
        \,=\,
        \sum_{k=1}^{K}
        \Bigg(
            \sigma_{k,0}
            + \sum_{i \in \mathcal{G}_k} g_i\sigma_{k,i}
            + \sum_{j \in \mathcal{H}_k} h_j\lambda_{k,j}
        \Bigg)
    \end{equation*}
\end{theorem}
With $p=f-\beta$, any strict lower bound can be certified by a sum of polynomials that are both nonnegative on $\mathcal{F}$ and \textit{local}---involving variables from just one clique each.
Consequently, no polynomial basis needs cross-clique monomials, e.g.\ $x_1x_3$ for \cref{eq:csp_example}.

This localization can be applied directly to our sample-based construction.
Generate samples $\mathcal{Z} \subset \mathcal{W}$ and extract $\mathcal{Z}_{C_k} := \{x_{C_k} \mid x \in \mathcal{Z}\}$ for each clique.
Form~$E_d(\mathcal{Z}_{C_k})$ using just the $\varpi(|C_k|,d)$ monomials in $x_{C_k}$, and compute a local quotient basis $\varphi_k : \mathbb{R}^{|C_k|} \to \mathbb{R}^{r_k}$ from its SVD as in \cref{sec:sample_based_sos}.
These bases can parameterize the local polynomials appearing in the certificate family of \cref{thm:spositiv}.
Evaluating the resulting identity on $\mathcal{Z}$ then gives the sparse analogue of \cref{eq:ineq_sdp},
\begin{equation}
\begin{gathered}
    \max \ \ \beta
    \quad
    \text{over } \, \beta \in \mathbb{R} \,
    \text{ and, for each } k \in \{1,\ldots,K\},
    \\[5pt]
    G_{k,i} \in \mathbb{S}^{r_k}_+\ \ \big(i \in \{0\} \cup \mathcal{G}_k\big)
    \qquad
    H_{k,j} \in \mathbb{S}^{r_k}\ \ \big(j \in \mathcal{H}_k \cap J\big)
    \\[5pt]
    \text{subject to the following at each sample }x \in \mathcal{Z}
    \\[3pt]
    f(x) - \beta
    \,=\,
    \sum_{k=1}^{K}
    \varphi_k(x_{C_k})^\top
    \Bigg(
        G_{k,0}
        + \sum_{i \in \mathcal{G}_k} g_i(x) G_{k,i}
        \,+\, \smashoperator{\sum_{j \in \mathcal{H}_k \cap J}} h_j(x) H_{k,j}
    \Bigg)
    \varphi_k(x_{C_k})
\end{gathered}
\label{eq:mono_sparse_sdp}
\end{equation}
Under the hypotheses of \cref{thm:spositiv}, increasing $d$ while maintaining the $\kappa$-poisedness condition from \cref{sec:handling_inequalities} gives bounds that converge monotonically to $f^*$.

By replacing the large PSD matrices in \cref{eq:ineq_sdp} with many smaller clique-sized blocks, \cref{eq:mono_sparse_sdp} already offers a substantial computational advantage.
However, this basic application of correlative sparsity leaves a critical bottleneck in the linear system itself.
Whereas traditional coefficient matching is inherently localized by monomial support, matching at a sample folds every clique into the same equation, even coupling blocks whose cliques share no variables.
So while the PSD cone is decomposed across cliques, the sampled linear system remains globally coupled.
Moreover, \cref{eq:mono_sparse_sdp} assumes that $\mathcal{Z}$ is $\kappa$-poised for $\mathcal{W}$ as a whole, demanding as many as $\varpi(n,\kappa)$ samples even when every clique is small.
To eliminate these bottlenecks, we introduce auxiliary polynomials that separate the sampled system into local equations coupled only where their cliques overlap.
We then establish when the resulting hierarchy preserves the bounds of \cref{eq:mono_sparse_sdp}.


\subsection{Separator Polynomials}
\label{sec:separator_polynomials}

Write $\beta = \sum_{k=1}^{K} \beta_k$ with $\beta_{1:K} \in \mathbb{R}$ and define the \textit{local residuals} by
\begin{equation*}
    \zeta_k
    \,:=\,
    f_k - \beta_k
    -
    \varphi_k^\top
    \Bigg(
        G_{k,0}
        + \sum_{i \in \mathcal{G}_k} g_i G_{k,i}
        \,+\, \smashoperator{\sum_{j \in \mathcal{H}_k \cap J}} h_j H_{k,j}
    \Bigg)
    \varphi_k
    \:\in\, \mathbb{R}[x_{C_k}]
\end{equation*}
Then the certificate system in \cref{eq:mono_sparse_sdp} is simply $\sum_{k=1}^{K} \zeta_k = 0$ on $\mathcal{Z}$.
We could decouple this system by setting each residual to zero individually, but doing so would exclude cancellations between nonzero residuals and likely weaken the bound.
Instead, we parameterize such cancellations with polynomials we call \textit{separators}.

Any monomial appearing in both $\zeta_k$ and $\zeta_l$ can only involve variables from the overlap $C_k \cap C_l$.
Thus, for every pair $k \leq l$, we represent the part of $\zeta_k$ canceled by $\zeta_l$ with a separator $\tau_{kl} \in \mathbb{R}[x_{C_k \cap C_l}]$ and require $\tau_{lk} = -\tau_{kl}$ for its counterpart in $\zeta_l$.\footnote{
    It follows that $\tau_{kk} = 0$ for any $k$, and that $\tau_{kl}$ is constant if $C_k \cap C_l$ is empty.
}
We then require each residual to match the sum of its separators at the samples,
\begin{equation}
    \zeta_k(x_{C_k}) \,=\, \sum_{l=1}^{K} \tau_{kl}(x_{C_k \cap C_l})
    \quad \forall\, x \in \mathcal{Z},\ k \in \{1,\ldots,K\}
    \label{eq:separated_system}
\end{equation}
Each separator pair can be viewed as a single transfer, subtracting $\tau_{kl}$ from $\zeta_k$ and adding it to $\zeta_l$.\footnote{
    The separators can also be understood as \textit{cost shifting}, redistributing the objective among its components before each is certified locally.
    E.g.\ in \cref{eq:csp_example} where $C_2 \cap C_3 = \{3\}$, we can write $f = \big(-x_2x_3 - \tau_{2,3}(x_3)\big) + \big(x_4^2 - x_5 - \tau_{3,2}(x_3)\big)$.
}
These pairs cancel when summed across cliques, maintaining zero total residual while enabling nonzero local residuals.
The question is whether such pairwise transfers can represent all possible cancellations among the local residuals.

\begin{theorem}[Membership Localization]
    \label{thm:membership_localization}
    Let $\mathcal{A} \subseteq \mathbb{R}[x_{1:n}]$, $C_{1:K} \subseteq \{1,\ldots,n\}$, and $\kappa \in \mathbb{N}$.
    Writing $\mathbb{R}[x_{1:n}]_{\leq\kappa}$ for the degree-$\leq\!\kappa$ polynomials in $x_{1:n}$, suppose
    \begin{equation}
        \mathcal{A}
        \,\cap\,
        \Bigg(
            \sum_{k=1}^{K} \mathbb{R}[x_{C_k}]_{\leq\kappa}
        \Bigg)
        \,=\,
        \sum_{k=1}^{K}
        \Big(
            \mathcal{A}
            \,\cap\,
            \mathbb{R}[x_{C_k}]_{\leq\kappa}
        \Big)
        \label{eq:ideal_decomp}
    \end{equation}
    Then, for any family $\zeta_k \in \mathbb{R}[x_{C_k}]_{\leq\kappa}$, the sum $\sum_k \zeta_k$ belongs to $\mathcal{A}$ if and only if there exist separators $\tau_{kl} \in \mathbb{R}[x_{C_k \cap C_l}]_{\leq\kappa}$ with $\tau_{lk} = -\tau_{kl}$ satisfying $\zeta_k - \sum_l \tau_{kl} \in \mathcal{A}\ \ \forall k$.
\end{theorem}
\begin{proof}
    Suppose $\sum_k \zeta_k \in \mathcal{A}$.
    By \cref{eq:ideal_decomp}, there exist $u_k \in \mathcal{A} \cap \mathbb{R}[x_{C_k}]_{\leq\kappa}$ such that $p_k := \zeta_k - u_k$ satisfy $\sum_k p_k = 0$.
    To construct separators satisfying $\sum_l \tau_{kl} = p_k$, initialize them all to zero.
    For each monomial $x^\alpha$ appearing in the $p_k$, let $c_k$ be its coefficient in $p_k$ and choose a nonzero $c_l$.
    For every $k \ne l$, add $c_k x^\alpha$ to $\tau_{kl}$ and subtract it from $\tau_{lk}$.
    This contributes coefficient $c_k$ to the $k$th separator sum for $k \ne l$, while the $l$th receives $-\sum_{k \ne l} c_k = c_l$, since $\sum_k c_k = 0$.
    Thus the separators have the required scopes and degrees and satisfy $\zeta_k - \sum_l \tau_{kl} = u_k \in \mathcal{A}$.
    Conversely, separator antisymmetry and \cref{eq:ideal_decomp} give $\sum_k \zeta_k = \sum_k \big(\zeta_k - \sum_l \tau_{kl}\big) \in \mathcal{A}$.
\end{proof}

Taking $\mathcal{A} = \{0\}$, condition \cref{eq:ideal_decomp} holds automatically and \cref{thm:membership_localization} shows that any zero-sum family of local degree-$\leq\!\kappa$ polynomials admits separators of the same degree.
Applying this to the certificates of \cref{thm:spositiv} and evaluating on $\mathcal{Z} \subset \mathcal{W}$ recovers \cref{eq:separated_system}.
These sampled local systems hold throughout $\mathcal{W}$ if each individual $\mathcal{Z}_{C_k}$ is $\kappa$-poised for $\mathcal{W}_{C_k} := \{x_{C_k} \mid x \in \mathcal{W}\}$, meaning that every degree-$\leq\!\kappa$ polynomial vanishing on the local samples also vanishes throughout this projection.\footnote{
    Even this condition can be stricter than necessary, as we'll see in \cref{sec:selecting_samples}.
}
Notably, the samples no longer need to distinguish polynomials across all $n$ variables, only those within each clique.
Then, summing \cref{eq:separated_system} across cliques certifies $\beta = \sum_k \beta_k \leq f^*$, while \cref{thm:spositiv} guarantees convergence to $f^*$ as $d$ is increased.

Furthermore, \cref{thm:membership_localization} identifies when using \cref{eq:separated_system} with local poisedness yields the same fixed-degree bound as \cref{eq:mono_sparse_sdp} with global poisedness.
Taking $\mathcal{A} = \mathcal{I}(\mathcal{W})$, condition \cref{eq:ideal_decomp} allows all vanishing relations in the total residual to be resolved within individual cliques, leaving only like-terms for the separators to cancel across cliques.
This condition is guaranteed if each clique contains the entire scope of every $h_{j\notin J}$ it shares variables with (e.g.\ quaternion normalization in multi-pose estimation).

The current pairwise construction for \cref{eq:separated_system} introduces $\binom{K}{2}$ separators, but far fewer can suffice to encode the exact same residual cancellations.
To see this, suppose three cliques $C_{1:3}$ all support a monomial $x^\alpha$ and that we wish to cancel $x^\alpha$ in $\zeta_1$ against $-x^\alpha$ in $\zeta_3$.
Instead of using $\tau_{1,3}=x^\alpha$, we can route it through $\tau_{1,2}=\tau_{2,3}=x^\alpha$ with no net impact on $\zeta_2$ because $\tau_{2,1}+\tau_{2,3}=-x^\alpha+x^\alpha=0$.
More generally, define the \textit{occurrence set} $\Omega_\alpha$ as the set of cliques supporting $x^\alpha$.
All pairwise cancellations involving $x^\alpha$ can be reproduced along any spanning tree on $\Omega_\alpha$, using only $|\Omega_\alpha|-1$ edges instead of all $\binom{|\Omega_\alpha|}{2}$ possible pairs.
We therefore choose a spanning tree for each distinct occurrence set and retain one separator pair $\tau_{kl}$ for each edge $\{k,l\}$ in the union of those trees.
To keep the number of separators small, we favor trees that reuse edges already selected for other occurrence sets.

RIP has so far only been used to ensure convergence through \cref{thm:spositiv}; neither bound validity nor occurrence-set routing requires it.
However, RIP greatly simplifies the separators by letting all occurrence sets share a single clique tree.
To see why, consider a leaf clique $C_k$ with parent $C_{\pi(k)}$.
RIP guarantees that any monomial $x^\alpha$ supported by $C_k$ and another clique is also supported by $C_{\pi(k)}$, so its coefficient at the leaf can be passed to the parent through $\tau_{k,\pi(k)}$.
Removing the leaf and repeating this process passes the coefficients of $x^\alpha$ all the way to the last remaining clique in $\Omega_\alpha$, where they sum to zero.
Thus all monomial cancellations can be routed through the same $K-1$ separators $\tau_{k,\pi(k)} \in \mathbb{R}[x_{C_k \cap C_{\pi(k)}}]$.
This efficiency presents a design decision when RIP does not hold: use occurrence-set routing with the existing cliques or enlarge them until RIP holds.\footnote{
    A \textit{chordal completion} adds edges to the CSP graph until its maximal cliques satisfy RIP.
    When this requires little fill-in, it can be preferable to occurrence-set routing.
}

With the routing chosen, all that remains is to construct a quotient basis for each separator.
As usual, take $\kappa \geq \max_k \deg(\zeta_k)$ uniformly.
Assign each supported monomial to the edges of its occurrence-set's spanning tree, and let $\mathcal{E}$ be the union of those trees.
For each edge $e = \{k,l\} \in \mathcal{E}$, the assigned monomials involve only $x_{C_e}$, where $C_e := C_k \cap C_l$, and span the polynomial space needed for its separator.\footnote{
    With RIP, $\mathcal{E}$ is exactly the clique tree and each edge carries all degree-$\leq\!\kappa$ monomials in $x_{C_e}$.
}
Take the corresponding columns of $E_\kappa(\mathcal{Z}_{C_e})$, center them to remove constants that can be absorbed into $\beta_{1:K}$, and use their SVD to compute a quotient basis $\psi_e : \mathbb{R}^{|C_e|} \to \mathbb{R}^{r_e}$.
Finally, write $\tau_{kl} = \varepsilon_{k,e}\, \gamma_e^\top \psi_e$ with a shared coefficient vector $\gamma_e \in \mathbb{R}^{r_e}$ and opposite endpoint signs $\varepsilon_{k,e} = -\varepsilon_{l,e} \in \{-1,1\}$.
The resulting fully sparsified SDP is
\begin{equation}
\begin{aligned}
    \max\ \sum_{k=1}^{K} \beta_k \quad \text{over}\quad\ &
    \beta_{k} \in \mathbb{R} \quad\;
    G_{k,i} \in \mathbb{S}^{r_k}_+ \quad\;
    H_{k,j} \in \mathbb{S}^{r_k} \quad\;
    \gamma_e \in \mathbb{R}^{r_e}
    \\
    \text{for all}\quad\ &
    i \in \{0\} \cup \mathcal{G}_k \quad\;
    j \in \mathcal{H}_k \cap J \quad\;
    e \in \mathcal{E} \quad\;
    k \in \{1,\ldots,K\}
    \\[7pt]
    \text{subject to}\quad\ &
    \zeta_k(x_{C_k})
    \;=\;\,
    \smashoperator{\sum_{\{e \in \mathcal{E} \mid k \in e\}}}
    \;\varepsilon_{k,e}\, \gamma_e^\top \psi_e(x_{C_e})
    \quad\ \forall\, x \in \mathcal{Z}
\end{aligned}
\label{eq:sparse_sdp}
\end{equation}

Each clique's dual sample weights $y_{k,1:S}$ now define a local functional $\mathcal{L}_k[p] := \sum_{s=1}^{S} y_{k,s} p(z_{s,C_k})$.
In the dual, $\beta_k$ enforces $\mathcal{L}_k[1] = 1$ while $\gamma_e$ enforces $\mathcal{L}_k[\psi_e] = \mathcal{L}_l[\psi_e]$ for $e = \{k,l\}$.
Through occurrence-set routing, these constraints ensure agreement of shared degree-$\leq\!\kappa$ moments across cliques, allowing the local first moments to be stitched into a candidate $\bar{x}$ for the usual refinement and gap assessment.

In comparison to \cref{eq:mono_sparse_sdp}, \cref{eq:sparse_sdp} needs only clique-local poisedness and replaces globally coupled sample equations with sparsely coupled local systems.
For problems with many small cliques and narrow overlaps, the $\gamma_e$ add modest free-variable overhead in exchange for a block-sparse constraint matrix that can be assembled clique by clique and exploited by most SDP solvers, with potential for distributed optimization.


\section{Selecting Samples}
\label{sec:selecting_samples}

As established in \cref{sec:leveraging_sparsity}, the poisedness requirement for \cref{eq:sparse_sdp} applies separately to each local sample set $\mathcal{Z}_{C_k} \subset \mathcal{W}_{C_k}$.\footnote{
    Each $\mathcal{Z}_{C_k}$ could even be generated independently, so long as the separator bases are chosen consistently across overlapping cliques.
    It is simply convenient to project a common $\mathcal{Z} \subset \mathcal{W}$.
}
This locality, combined with the freedom to choose $\mathcal{W}$ through $J$, makes poisedness much easier to achieve via random sampling.
The remaining task is to determine just how many samples are actually needed for a given problem.

Let $\mathcal{R}_k \subseteq \mathbb{R}[x_{C_k}]_{\leq\kappa}$ be the span of all possible \textit{adjusted residuals} $\zeta_k - \sum_l \tau_{kl}$.
Local poisedness amounts to every polynomial $p \in \mathcal{R}_k$ being uniquely determined on $\mathcal{W}$ by its sampled values, so that $p \in \mathcal{I}(\mathcal{Z}_{C_k}) \,\Rightarrow\, p \in \mathcal{I}(\mathcal{W})$.
These polynomials are specified by $D_k := \dim\big(\mathcal{R}_k|_{\mathcal{W}}\big)$ coefficients, while each sample supplies one linear equation in those coefficients.
Thus at least $D_k$ samples are needed, and they suffice precisely when the sampled system has rank $\rho_k = D_k$.
One could find $D_k$ algebraically, e.g.\ via Gr\"obner bases, but that is often intractable.
Instead, we increase the sample count until linear redundancy reveals that each $\rho_k$ has reached $D_k$, then prune the excess.

Let $\Phi_{k,s}$ denote row $s$ of the local product evaluation matrix $\Phi_k(\mathcal{Z}_{C_k})$, defined analogously to $\Phi(\mathcal{Z})$ in \cref{sec:sample_based_sos}.
Define row $s$ of the $S \times \Theta_k$ matrix $A_k$ as
\begin{equation*}
    [A_k]_{s,:} \,:=\, \Big[\, 1 \,\textcolor{gray!75}{\big\vert}\, \Phi_{k,s} \,\textcolor{gray!75}{\big\vert}\, \big(g_i(z_s)\Phi_{k,s}\big)_{i \in \mathcal{G}_k} \,\textcolor{gray!75}{\big\vert}\, \big(h_j(z_s)\Phi_{k,s}\big)_{j \in \mathcal{H}_k \cap J} \,\textcolor{gray!75}{\big\vert}\, \big(\varepsilon_{k,e}\psi_e(z_{s,C_e})^\top\big)_{e \in \mathcal{E} \,:\, k \in e} \,\Big]
\end{equation*}
Then the sampled system for each clique $k$ reads $A_k\theta_k = b_k$, where $\theta_k \in \mathbb{R}^{\Theta_k}$ vectorizes the corresponding SDP decision variables and $b_k = f_k(\mathcal{Z}_{C_k})$.
As $\theta_k$ varies, the adjusted-residual evaluation vectors $b_k - A_k\theta_k$ have linear span equal to the column space of $\overline{A}_k := [A_k\ b_k]$.
Hence $\rho_k = \operatorname{rank}\overline{A}_k$ can be computed numerically.

To determine when $\rho_k$ has reached $D_k$, assume $J$ has been chosen so that each irreducible component of $\mathcal{W}$ can be sampled readily.
On each component, samples drawn independently from a nondegenerate\footnote{
    I.e., the distribution assigns zero probability to any proper algebraic subset of that component.
}
distribution almost surely increase the rank of the corresponding row block of each $\overline{A}_k$ by one per draw until the adjusted residuals are fully determined there, as in \mbox{\cite[Section~5]{cifuentes2017sampling}}.
We therefore draw samples from each component until its row block in every $\overline{A}_k$ has more rows than its rank, ensuring that the combined samples satisfy $\rho_k = D_k$ for all $k$.
If the chosen $\mathcal{W}$ is irreducible, this saturation procedure is simply to increase $S$ until each $\rho_k < S$.

The largest clique size $w := \max_k |C_k|$ gives $\max_k D_k \leq \varpi(w,\kappa)$, but rather than start with that count, we grow the sample set in modest batches, repeating the cheap rank checks each time.
Once all ranks saturate, we use a pivoted QR of $\smash{\overline{A}_k^{\raisebox{-0.5ex}{$\scriptstyle\top$}}}$ to select $\rho_k = D_k$ independent rows of $\overline{A}_k$.
With $P_k \in \{0,1\}^{D_k \times S}$ as the selector,
\begin{equation*}
    \textcolor{gray}{\underbrace{\normalcolor A_k\theta_k = b_k}_{S\text{ equations}}}
    \quad\iff\quad
    \textcolor{gray}{\underbrace{\normalcolor (P_k A_k)\theta_k = P_k b_k}_{\mathclap{D_k \leq S\text{ equations, same solutions}}}}
\end{equation*}
Retaining original equations (rather than, say, an SVD row-space basis) preserves the rank-one constraint-matrix structure\footnote{
    Not to be confused with the rank-one optimal moment matrices used for solution extraction.
}
$\varphi_k(z_{s,C_k}) \varphi_k(z_{s,C_k})^\top$ exploited by some SDP solvers \cite{lofberg2004sampling} and keeps the dual variables directly interpretable as sample weights.

Returning to \cref{eq:doublecircle}, draw $100$ samples from each circle in $\mathcal{V}$ and let $d=3,\,K=1$.
The $12$-element quotient basis gives $\overline{A} = [\mathbf{1}_{200}\ \Phi(\mathcal{Z})\ f(\mathcal{Z})]$ with $80$ columns, but its rank saturates at just $D = 24$.
Pivoted QR retains $24$ sample equations, recovering the $24\times79$ system in \cref{sec:sample_based_sos}.
Using samples from the ambient variety of \cref{sec:circumventing_enumeration} instead, QR retains $73$ equations, revealing that $\mathcal{R}|_{\mathcal{W}}$ has dimension $D = 73 > 24$.

Note that while the SDP is constructed from random samples, poisedness ensures that its optimal solution set is deterministically fixed up to a known change of coordinates.
I.e., sample-based SOS is not a stochastic optimizer.
That said, the choice of sampling distribution can affect numerical conditioning in finite-precision arithmetic \cite{lofberg2004sampling}.
As a practical heuristic, we suggest centering and scaling the distribution around the best available estimates of the global optima, while deliberately retaining some spread to distinguish the relevant polynomials reliably.


\section{Results}
\label{sec:results}

We provide \textit{Rhiza}, a high-performance Python implementation of our complete methodology:
\textcolor{blue}{\href{https://github.com/Rhiza-PolyOpt}{\texttt{https://github.com/Rhiza-PolyOpt}}}\\
Specifying a problem requires only the variable scope and sampled data of each objective component and constraint.
Rhiza automatically constructs and prunes the sparse SOS relaxation \cref{eq:sparse_sdp}, passes it to an off-the-shelf SDP solver, and extracts solution candidates, with careful numerical conditioning and diagnostics integrated throughout.
Its fully sample-based interface also supports a raw data format, allowing users to generate their problem from models in any programming language.
A detailed user manual and test suite support its use in new applications.

We demonstrate the practical reach of this implementation on realistic instances of 3D point-cloud registration and spacecraft thruster scheduling, two challenging engineering problems that combine nonlinear geometry with many discrete decisions.
At scales involving thousands of continuous and discrete variables, we obtain solutions in seconds with certifiably small optimality gaps throughout.

Both applications use polynomial lifts to obtain more expressive certificates at low $d$.
In sample-based SOS, a lift $p$ simply extends each sample $z_s$ to $(z_s,p(z_s))$, replacing the unreduced basis $\nu_d(x)$ with $\nu_d(x,p(x))$.
At $d=1$, for instance, quadratic lifts provide selected degree-$2$ polynomials without incurring the entire $d=2$ basis in $x$.
Sample-based quotienting then makes full use of this enrichment without the symbolic substitutions and redundant constraints generally needed by traditional SOS.

All reported experiments were performed on a Dell Precision laptop with an Intel Core i7-12800H processor and 31\,GiB of usable RAM.
We used the standard SDP solver \href{https://www.mosek.com/}{MOSEK} with 14 threads and all default settings except for the stopping tolerance, which was chosen once for each problem family based on a single instance.\footnote{
    The experiment code and configurations are provided at the same link.
    Rhiza also supports the open-source SDP solver \href{https://www.cvxgrp.org/scs/}{SCS}, so a MOSEK license is not required.
}


\subsection{Geometric Perception}
\label{sec:geometric_perception}

Point-cloud registration (PCR) is a foundational problem in geometric perception, tightly coupling discrete data association with continuous pose estimation.
Given two unordered point clouds $\{a_i\}_{i=1}^{P},\,\{b_j\}_{j=1}^{P} \subset \mathbb{R}^3$ related by an unknown rigid pose and Gaussian measurement noise, we seek the rotation $R \in \mathrm{SO}(3)$, translation $t \in \mathbb{R}^3$, and correspondences between their points.
Since these clouds represent scene geometry rather than labeled landmarks, each point may match any points in the other cloud, motivating an ``all-to-all'' truncated least-squares model with binary variables $w_{ij}$ for all $P^2$ possible correspondences.
\vspace{-0.25em}
\begin{equation}
\begin{aligned}
    \min_{R,t,w}\quad
    & \frac{1}{P^2}\sum_{i,j=1}^{P}\Big(w_{ij}\|Ra_i + t - b_j\|^2 + (1-w_{ij})c^2\Big) + \frac{1}{P}\|R\bar{a} + t - \bar{b}\|^2
    \\[5pt]
    \text{s.t.}\quad
    & R^\top\!R = I
    \,,\quad
    \det(R) = 1
    \,,\quad
    \|t\|^2 \leq D^2
    \,,\quad
    w_{ij}(1 - w_{ij}) = 0 \;\ \forall\,i,j
\end{aligned}
\label{eq:pcr}
\end{equation}
Each $w_{ij}\!\in\!\{0,\!1\}$ selects between a Gaussian-inlier squared error and a uniform-outlier rejection penalty $c^2$.
The centroid regularizer and translation bound provide weak priors, encouraging mean alignment and limiting translation to the scene diameter $D$.

We form $K=P^2$ fixed-size cliques: one per binary $w_{ij}$, all sharing the pose $(R,t)$ and arranged in a star tree satisfying RIP.
We lift the $d=1$ local bases with the \mbox{degree-$2$} products $w_{ij}R$ and $w_{ij}t$ by appending their values to each sample.
Although $\mathcal{V}$ has $2^{\scriptscriptstyle{P^2}}$ components, each clique sees only two, so rank saturation shows that just $142$ samples suffice regardless of $K$.\footnote{
    All PCR experiments here use the same $160 > 142$ samples for numerical margin.
}
Pose rounding and Hungarian assignment recover a solution from the first moments.
Performance is reported in \cref{tab:pcr_results}.

\begin{table}[t]
    \centering
    \caption{
    PCR results with $K=P^2$ candidate correspondences, only $P$ of which are true (an inlier rate of $1/P$).
    Despite inlier rates as low as $2.5\%$, first moments indicate over $95\%$ of the true matches.
    Note that measurement noise generally shifts the optimal pose away from the ground truth.
    }
    \label{tab:pcr_results}
    \small
    \setlength{\tabcolsep}{4pt}
    \begin{tabular}{@{}cccccc@{}}
        \toprule
        \shortstack{Candidate\\Correspondences} &
        \shortstack{Build\\Time} &
        \shortstack{Solve\\Time} &
        \shortstack{Relative Optimality\\Gap $\;\Delta_\varphi(\hat{x})/D^2$} &
        \shortstack{Relative Translation\\Error $\;\|\hat{t}-t_{\scriptscriptstyle{\mathrm{truth}}}\|/D$} &
        \shortstack{Rotation\\Error} \\
        \midrule
        $10^2 = \hphantom{1}100$ & $0.53\,\mathrm{s}$ & $2.6\,\mathrm{s}$ & $0.0060\%$ & $0.67\%$ & $1.0^\circ$ \\
        $20^2 = \hphantom{1}400$ & $1.9\,\mathrm{s}$  & $11\,\mathrm{s}$  & $0.041\%$  & $0.61\%$ & $4.3^\circ$ \\
        $30^2 = \hphantom{1}900$ & $4.3\,\mathrm{s}$  & $22\,\mathrm{s}$  & $0.057\%$  & $0.30\%$ & $4.5^\circ$ \\
        $40^2 = 1600$            & $8.0\,\mathrm{s}$  & $40\,\mathrm{s}$  & $0.073\%$  & $0.36\%$ & $6.7^\circ$ \\
        \bottomrule
    \end{tabular}
    \vspace{-1em}
\end{table}

Our build and solve times \textbf{grow linearly} in $K$, reflecting the fixed-size cliques and enabling us to handle $1600$ correspondences in under $49$ seconds.
For comparison, Yang and Carlone \cite[Table~II]{yang2022outliers} report an average MOSEK solve time of $377$ seconds for their hand-crafted SDP at just $20$ correspondences, with memory exhaustion at $100$.
While they manually designed their basis and derived redundant constraints $h_j$, Rhiza automatically built our SDP from simple pose and binary samples, inherently capturing the algebraic consequences of every possible $h_j$ at all degrees.


\subsection{Motion Planning}
\label{sec:motion_planning}

A spacecraft must reorient through five timed attitude waypoints while keeping angular speed below $1\,\mathrm{rad/s}$ and conserving fuel.
Its six thruster torque couples are either off or throttled between $25\%$ and $100\%$, coupling discrete firing decisions with continuous command magnitudes.\footnote{
    This deadband approximates minimum on-time effects in pulsed cold-gas thrusters.
}
Discretizing the $17.5\,\mathrm{s}$ maneuver into $100$ steps yields $3612$ lifted variables, including $600$ binaries, with state-transition cliques forming a chain satisfying RIP.
Just $360$ random-actuation rollouts encode the dynamics constraints for the lifted $d=1$ relaxation, while mission and operating limits remain explicit inequalities.\footnote{
    The complete formulation and visuals are given in the experiments codebase documentation.
}

Rhiza builds the SDP in $4.5$ seconds and MOSEK solves it in $41$ seconds.
Local refinement yields a feasible trajectory in $3.3$ seconds with cost $114$, just $\Delta_\varphi(\hat{x}) \approx 3$ above the certified lower bound of $111$.
Meanwhile, a conventional boost--coast--brake heuristic with the same refinement yields cost $147$ and still slightly violates the speed limit.
For context, \cite[Table~1]{kang2024trajopt} reports $1648$ seconds to solve a $60$-step 2D flying-robot problem using a custom GPU-accelerated SDP solver, while Rhiza returns our $100$-step mixed-integer 3D maneuver in under $50$ seconds using MOSEK on a CPU.
Both exploit correlative sparsity, but Rhiza also quotients out the dynamics through sampling, enabling smaller SDPs with tighter bounds at lower degree.


\section{Conclusion}
\label{sec:conclusion}

As with other relaxation hierarchies, tightness is not guaranteed at tractable degrees.
Nevertheless, \cref{sec:results} shows how our novel combination of sampling and sparsity brings large mixed-integer semialgebraic problems within reach of a standard SDP solver.
Local quotient bases and separator polynomials tie SDP block structure and sampling requirements to clique complexity rather than total problem size.
Combining this formulation with GPU-accelerated SDP solvers and/or low-rank Burer--Monteiro methods could extend these gains to even larger problems and higher degrees.




\section*{Acknowledgments}

This work was funded by MIT Lincoln Laboratory and the OUSD(R\&E) through Air Force Contract FA8702-25-D-B002.
ChatGPT 6.0 Astra was used to revise the manuscript and help develop the codebase.

\bibliographystyle{siamplain}
\bibliography{references}

\end{document}